\documentclass{article}

\pdfoutput=1

\usepackage[utf8]{inputenc}
\usepackage{supertabular}
\usepackage{diagbox}%
\usepackage{float}
\usepackage{epstopdf}
\usepackage{caption}
\usepackage{titlesec}%
\usepackage{subfig}
\usepackage{graphicx}
\usepackage{indentfirst}
\usepackage{url}
\usepackage{listings}
\usepackage{appendix}
\usepackage{amsmath}
\usepackage{amssymb}
\usepackage{makecell}
\usepackage{multirow}
\usepackage{lipsum}
\usepackage{cite}

\usepackage{amsthm}
\usepackage{tabularx}

\usepackage{array}
\usepackage{booktabs}
\usepackage{caption}

\newtheorem{theorem}{Theorem}

\newtheorem{remark}{Remark}

\usepackage[T1]{fontenc}
\usepackage[utf8]{inputenc}
\usepackage{authblk}
\usepackage{geometry}
\newcommand\keywords[1]{\textbf{Keywords}: #1}

\title{An Efficient Multiple Harmonic Balance Method for Computing Quasi-Periodic Responses of Nonlinear Systems}
\author[a, b]{Qisi Wang}
\author[a, b]{Zipu Yan}
\author[a, b]{Honghua Dai \thanks{Corresponding author: hhdai@nwpu.edu.cn}}

\affil[a]{School of Astronautics, Northwestern Polytechnical University, 127 West Youyi Road, Xi'an, 710072, Shaanxi, China}
\affil[b]{National Key Laboratory of Aerospace Flight Dynamics, Northwestern Polytechnical University, 127 West Youyi Road, Xi'an, 710072, Shaanxi, China}

\date{\today}

\begin{document}

\maketitle

\begin{abstract}
Quasi-periodic responses composed of multiple base frequencies widely exist in science and engineering problems. The multiple harmonic balance (MHB) method is one of the most commonly used approaches for such problems. However, it is limited by low-order estimations due to complex symbolic operations in practical uses. Many variants have been developed to improve the MHB method, among which the time domain MHB-like methods are regarded as crucial improvements because of their high efficiency and simple derivation. But there is still one main drawback remaining to be addressed. The time domain MHB-like methods negatively suffer from non-physical solutions, which have been shown to be caused by aliasing (mixtures of the high-order into the low-order harmonics). Inspired by the collocation-based harmonic balancing framework recently established by our group, we herein propose a reconstruction multiple harmonic balance (RMHB) method to reconstruct the conventional MHB method using discrete time domain collocations. Our study shows that the relation between the MHB and time domain MHB-like methods is determined by an aliasing matrix,  which is non-zero when aliasing occurs. On this basis, a conditional equivalence is established to form the RMHB method. Three numerical examples demonstrate that this new method is more robust and efficient than the state-of-the-art methods.
\end{abstract}

\keywords{reconstruction multiple harmonic balance method,de- aliasing, sampling rules, quasi-periodic responses. }

\section{Introduction}
Quasi-periodic responses consist of multiple harmonic terms, which are not necessarily integer multiples of each other, commonly occurring in many scientific \cite{Bohr2018,Broer2009,Martin2010,Fan2016} and engineering studies \cite{Kim1998,Ushida1984,Pušenjak2004,Li2021,Chen2022}. Numerous ready-to-use methods have been developed for effectively solving such problems in the past few decades. Numerical integration methods are commonly adopted but limited by undesired simulation time for transient motion, small step size to constrain accumulated error and incapability of obtaining unstable periodic solutions. The parameter expansion methods (such as the perturbation method \cite{Bussgang1974} and averaging method \cite{Pengcheng1999,Shen2013}) are also available tools, but they presume the system to be weakly nonlinear, leading to limited applications. Thus, the semi-analytical methods, represented by the harmonic balance (HB) method \cite{Tseng1970} and its improved methods (single base frequency is selected) \cite{Lau1981,Cameron1989,Hall2002,Krack2019,Dai2022} are accepted as ideal options for solving nonlinear periodic response. However, the HB-like methods are time-consuming when computing complex periodic responses, because excessive harmonics have to be adopted to span the spectrum containing the low and high frequencies. What is more, there is no strict period for some quasi-periodic responses, so the presumed base frequency should be small for capturing a relatively accurate solution \cite{Chua1981}.

More efforts have been made to improve the efficiency of HB-like methods. The intuitive attempt to overcome the above restrictions is the MHB method. It supplements more base harmonics and works by presuming a multidimensional Fourier series for the desired periodic solutions and then computing resultant nonlinear algebraic equations (NAEs) of the frequency coefficients through balancing harmonics up to the truncation order \cite{Liu2007}. However, the derivation of MHB algebraic systems is still tedious and lengthy due to multidimensional Fourier analysis, hence a series of improved MHB-like were proposed \cite{Chua1981,Lau1983,Kim1996,Liu2007,Kundert1988}. Lau \cite{Lau1983} introduced the concept of multiple time scales to extend the application of the incremental harmonic balance (IHB) method. Although the IHB method can slightly reduce symbolic operations, the higher computation complexity in finding the Jacobian matrix and the residual vector using multiple integrals imposes limitations on its applications. Applying the two-dimensional discrete Fourier transform (DFT) to replace those integrals into trigonometric summations, Kim and Noah \cite{Kim1996,Kim1997} proposed the multiple harmonic balance alternating frequency/time (MHB-AFT) method. The MHB-AFT and the two-time-scales IHB method are shown to be mathematically equivalent \cite{Ju2020}, so the sophisticated derivation of the explicit Jacobian matrix and nonlinear force vector are still unavoidable. Since many trigonometric summations need to be calculated as the order of the MHB-AFT method increases, the computation becomes cumbersome in solving problems with complex nonlinearities. In order to make the computation more concise, Liu \cite{Liu2007} applied the high dimensional harmonic balance (HDHB) method with two base frequencies to compute quasi-periodic motions. The HDHB method is essentially an easy-to-implement time domain MHB-like method \cite{Dai2012}, that is to replace the original MHB algebraic system with simple time domain quantities. But its accuracy is impaired by aliasing-induced non-physical solutions \cite{Liu2006,Dai2014}. Many studies have been devoted to de-aliasing techniques, which can be divided into two main categories: numerical filtering de-aliasing \cite{Orszag1971,Labryer2009,Huang2014} and collocation-manipulating de-aliasing \cite{Shannon1949,Ju2020,Dai2022}. Recent studies have shown that the RHB methods can completely eliminate the aliasing of polynomial and rational fraction nonlinearities \cite{Dai2022}. Although the collocation-based RHB method is highly efficient, it is originally designed for periodic responses. As for how to avoid non-physical solutions in multiple base frequencies, there is still no systematic conclusion. 

Inspired by the collocation-based harmonic balancing framework established by our group \cite{Dai2022}, this study aims to propose a computationally-cheap reconstruction multiple harmonic balance (RMHB) method. In this study, the mechanism of aliasing in multiple harmonic balance computations is revealed and expressed by one aliasing matrix. It further leads to the discovery of a conditional identity,  which bridges the gap between the frequency domain and time domain multidimensional Fourier analysis. Here we show that the RMB method can be equivalently transformed into the original MHB method by choosing corresponding sampling rules when computing different kinds of quasi-periodic responses. The present study successfully solves some computational difficulties of the HB-like methods (low computing efficiency due to excessive harmonics), the MHB method (lengthy symbolic operations), and the MHB-like time domain methods (aliasing phenomenon and oversampling problem).

Because the dynamic response usually contains two base harmonics in engineering studies, we mainly discuss two base frequencies situations in this paper. The performance of the RMHB method is evaluated by three nonlinear examples from structural vibrations, microwave circuits, to aeroelasticity dynamics problem. First, the features of the quasi-periodic computation are explored using the classical forced Van der Pol equation. Then the Duffing oscillator with multiple input frequencies is investigated using the RMHB method to demonstrate both its efficiency and accuracy in dealing with complex periodic responses without causing aliasing. Moreover, the nonlinear aeroelastic system of an airfoil with an external store is studied. The complicated responses will be obtained semi-analytically by employing the RMHB method. Finally, the numerical results verified the effectiveness of the proposed RMHB method. Although the number of harmonic components grows quadratically with increasing truncation order, this is a common feature of all MHB-like methods. This new method still guarantees a good trade-off between efficiency and accuracy in contrast with the existing methods.

\section{Reconstruction Multiple Harmonic Balance Method}\label{}
The ordinary differential equations for a general $N$-DOF nonlinear dynamical system can be expressed as
\begin{equation}
    \dot{\mathbf{x}}=\mathbf{f}(\mathbf{x},t),\label{eq:original}
\end{equation}
where the state vector $\mathbf{x}=[ x_1,\cdots ,x_N ] ^{\mathrm{T}}$. For the two base frequencies cases, every single DOF can be written in the following form \cite{Liu2007}:
\begin{equation}
    x_{i}(t)=\sum_{m} \sum_{n} \hat{x}_{i c}(m, n) \cos \left(\left(m \omega_{1}+n \omega_{2}\right) t\right) +\hat{x}_{i s}(m, n) \sin \left(\left(m \omega_{1}+n \omega_{2}\right) t\right),
\end{equation}
and parameters $m$ and $n$ satisfy 
\begin{equation}
    \left| m \right|+\left| n \right|\leqslant p,
    \label{relation}
\end{equation}
where $p$ is the truncation order for 2D Fourier series \cite{Ju2017}. Thus the state vector and nonlinear terms can be approximated as
\begin{equation}
\begin{aligned}
    &\mathbf{x}(t)=\mathbf{I}_N\otimes \left[ 1\,\,\mathrm{c}^{1,0}\,\,\mathrm{s}^{1,0}\cdots \mathrm{c}^{m,n}\,\,\mathrm{s}^{m,n}\cdots \mathrm{c}^{0,p}\,\,\mathrm{s}^{0,p} \right] \mathbf{\hat{x}},\\
     &\mathbf{f}(\mathbf{x},t)=\mathbf{I}_N\otimes \left[ 1\,\,\mathrm{c}^{1,0}\,\,\mathrm{s}^{1,0}\cdots \mathrm{c}^{m,n}\,\,\mathrm{s}^{m,n}\cdots \mathrm{c}^{0,p}\,\,\mathrm{s}^{0,p} \right] \mathbf{\hat{f}}+O(\mathbf{\hat{f}}) .\label{eq:represents}      
\end{aligned}
\end{equation}
The vector of unknowns is
\begin{equation}
    \mathbf{\hat{x}}=[ \hat{x}_{1c}(0,0) ,\cdots ,\hat{x}_{ic}(m,n) ,\hat{x}_{is}(m,n) ,\cdots ,\hat{x}_{Nc}(0,p) ,\hat{x}_{Ns}(0,p)] ^{\mathrm{T}},\nonumber
\end{equation}
where $\mathbf{\hat{f}}$ is a nonlinear polynomial function of $\mathbf{\hat{x}}$. $\mathbf{I}_N$ is the identity matrix of dimension $N$, and $\otimes$ is the Kronecker product. Denote $\mathrm{c}^{m,n}=\cos ( m\omega _1+n\omega _2) t$ and $\mathrm{s}^{m,n}=\sin(m\omega _1+n\omega _2) t$. Besides the first derivative of $\mathbf{x}(t)$ with respect to time is expressed \cite{Krack2019}
\begin{equation}
    \mathbf{\dot{x}}(t)=\mathbf{I}_N\otimes \left[ 1\,\,\mathrm{c}^{1,0}\,\,\mathrm{s}^{1,0}\cdots \,\,\mathrm{c}^{m,n}\,\,\mathrm{s}^{m,n}\cdots \,\,\mathrm{c}^{0,p}\,\,\mathrm{s}^{0,p} \right] \nabla \mathbf{\hat{x}},\label{eq:derivative}
\end{equation}
with
\begin{equation}
    \nabla =\mathbf{I}_N\otimes \mathrm{diag}\left[ 0,\nabla _{1,0},\cdots ,\nabla _{m,n},\cdots ,\nabla _{0,p} \right] ,\,\,
    \nabla _{m,n}=\left[ \begin{matrix}
	0&		m\omega _1+n\omega _2\\
	-(m\omega _1+n\omega _2)&		0\\
\end{matrix} \right]. \nonumber
\end{equation}

By substituting Eqs. (\ref{eq:represents}) and (\ref{eq:derivative}) into (\ref{eq:original}), and applying the Galerkin method yields $N\times[2p(p+1)+1]$ NAEs of the MHB method
\begin{equation}
    \nabla \mathbf{\hat{x}}=\mathbf{\hat{f}}\left( \mathbf{\hat{x}} \right)\label{MHB NAEs},
\end{equation}
which can be readily solved by a NAE solver. However, it will get complicated and tedious for carrying the formula about nonlinear terms $ \mathbf{\hat{f}}$ as the truncation order of the MHB method increases.

Similar to the RHB method \cite{Dai2022}, the RMHB method establishes the relation between the Fourier coefficients and temporal quantities at $M$ equally spaced nodes over one period $T$ via a constant collocation matrix.
\begin{equation}
    \mathbf{\tilde{x}}_{M}=\mathbf{E\hat{x}},
\end{equation}
where $\mathbf{\tilde{x}}_{M}=\left[ \mathbf{x}\left( t_1 \right) \mathbf{x}\left( t_2 \right) \cdots \mathbf{x}\left( t_M \right) \right] ^{\mathrm{T}}$, $t_i=T(i-1)/M,i=1,2,\cdots,M$. The collocation matrix $\mathbf{E}$ and transformation matrix $\mathbf{E}^*$ are
\begin{equation}
\,\begin{cases}
	\mathbf{E}=\mathbf{I}_N\otimes \left[ \begin{matrix}
	1&		1&		\cdots&		1\\
	\mathrm{c}^{1,0}\left( t_1 \right)&		\mathrm{c}^{1,0}\left( t_2 \right)&		\cdots&		\mathrm{c}^{1,0}\left( t_M \right)\\
	\mathrm{s}^{1,0}\left( t_1 \right)&		\mathrm{s}^{1,0}\left( t_2 \right)&		\cdots&		\mathrm{s}^{1,0}\left( t_M \right)\\
	\vdots&		\vdots&		\cdots&		\vdots\\
	\mathrm{c}^{m,n}\left( t_1 \right)&		\mathrm{c}^{m,n}\left( t_2 \right)&		\cdots&		\mathrm{c}^{m,n}\left( t_M \right)\\
	\mathrm{s}^{m,n}\left( t_1 \right)&		\mathrm{s}^{m,n}\left( t_2 \right)&		\cdots&		\mathrm{s}^{m,n}\left( t_M \right)\\
	\vdots&		\vdots&		\cdots&		\vdots\\
	\mathrm{c}^{0,p}\left( t_1 \right)&		\mathrm{c}^{0,p}\left( t_2 \right)&		\cdots&		\mathrm{c}^{0,p}\left( t_M \right)\\
	\mathrm{s}^{0,p}\left( t_1 \right)&		\mathrm{s}^{0,p}\left( t_2 \right)&		\cdots&		\mathrm{s}^{0,p}\left( t_M \right)\\
\end{matrix} \right] ^{\mathrm{T}}\\
	\mathbf{E}^*=\mathbf{I}_N\otimes \frac{2}{M}\left[ \begin{matrix}
	\frac{1}{2}&		\frac{1}{2}&		\cdots&		\frac{1}{2}\\
	\mathrm{c}^{1,0}\left( t_1 \right)&		\mathrm{c}^{1,0}\left( t_2 \right)&		\cdots&		\mathrm{c}^{1,0}\left( t_M \right)\\
	\mathrm{s}^{1,0}\left( t_1 \right)&		\mathrm{s}^{1,0}\left( t_2 \right)&		\cdots&		\mathrm{s}^{1,0}\left( t_M \right)\\
	\vdots&		\vdots&		\cdots&		\vdots\\
	\mathrm{c}^{m,n}\left( t_1 \right)&		\mathrm{c}^{m,n}\left( t_2 \right)&		\cdots&		\mathrm{c}^{m,n}\left( t_M \right)\\
	\mathrm{s}^{m,n}\left( t_1 \right)&		\mathrm{s}^{m,n}\left( t_2 \right)&		\cdots&		\mathrm{s}^{m,n}\left( t_M \right)\\
	\vdots&		\vdots&		\cdots&		\vdots\\
	\mathrm{c}^{0,p}\left( t_1 \right)&		\mathrm{c}^{0,p}\left( t_2 \right)&		\cdots&		\mathrm{c}^{0,p}\left( t_M \right)\\
	\mathrm{s}^{0,p}\left( t_1 \right)&		\mathrm{s}^{0,p}\left( t_2 \right)&		\cdots&		\mathrm{s}^{0,p}\left( t_M \right)\\
\end{matrix} \right]\\
\end{cases}
\end{equation}

In the RMHB method, we take $\mathbf{\hat{f}}(\hat{x})=\mathbf{E}^*\mathbf{\tilde{f}}_M\left( \mathbf{\tilde{x}} \right)$ to replace the original nonlinear terms, where $\mathbf{\tilde{f}}_M\left( \mathbf{\tilde{x}} \right)$ is the value of $\mathbf{f}(\mathbf{x},t)$ at $M$ discrete time collocations $\tilde{\mathbf{x}}$. Here Eq. (\ref{MHB NAEs}) can be rewritten as
\begin{equation}
    \mathbf{E}^*\mathbf{E}\nabla \mathbf{\hat{x}}=\mathbf{E}^*\mathbf{\tilde{f}}_M\left( \mathbf{E\hat{x}} \right),\label{RMHB NAEs}
\end{equation}

If there are more than two base frequencies $\left\{ \omega _1,\omega _2,\cdots ,\omega _t \right\} $, the $p$-order RMHB method contains harmonic combinations as $m_{1p}\omega _1+m_{2p}\omega _2+\cdots +m_{tp}\omega _t$, and these coefficients of the linear combination satisfy the inequality relation \cite{Chua1981}:
\begin{equation}
    \left| m_{1p} \right|+\left| m_{2p} \right|+\cdots +\left| m_{tp} \right|\le p.
\end{equation}

\begin{figure}
	\centering
		\includegraphics[scale=0.25]{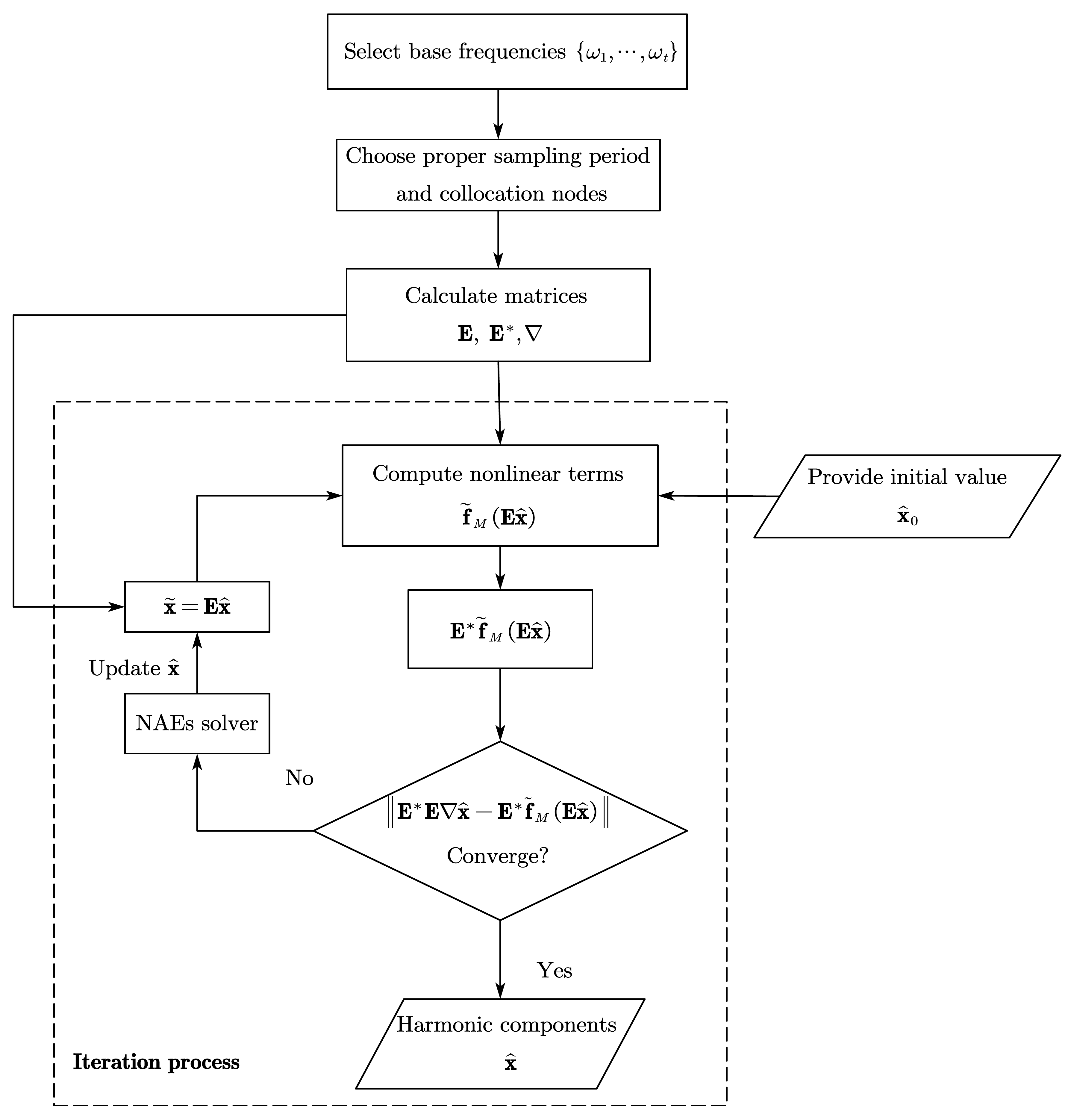}
		\caption{Computing procedure of the RMHB method.}
		\label{FIG:flow}
\end{figure}

Because quasi-periodic responses are a generalization of periodic responses, the single-base-frequency methods are seen to be special cases of the RMHB method. In order to better illustrate the method implementation, the flow chart of the RMHB method is sketched in Figure \ref{FIG:flow}. As for how to select the proper $T$ and $M$ to achieve a desirable accuracy, we will further explain this in the next section.
\section{Conditional Equivalence of RMHB and MHB}
The accuracy of the estimation solution will be affected by aliasing \cite{Dai2022,Dai2014}. To eliminate the aliasing in the multiple harmonic balance computations, here we conclude two sets of sampling rules for $T$ and $M$ for different multiple frequency harmonic balance computations.

\begin{figure}
	\centering
	 \includegraphics[scale=0.35]{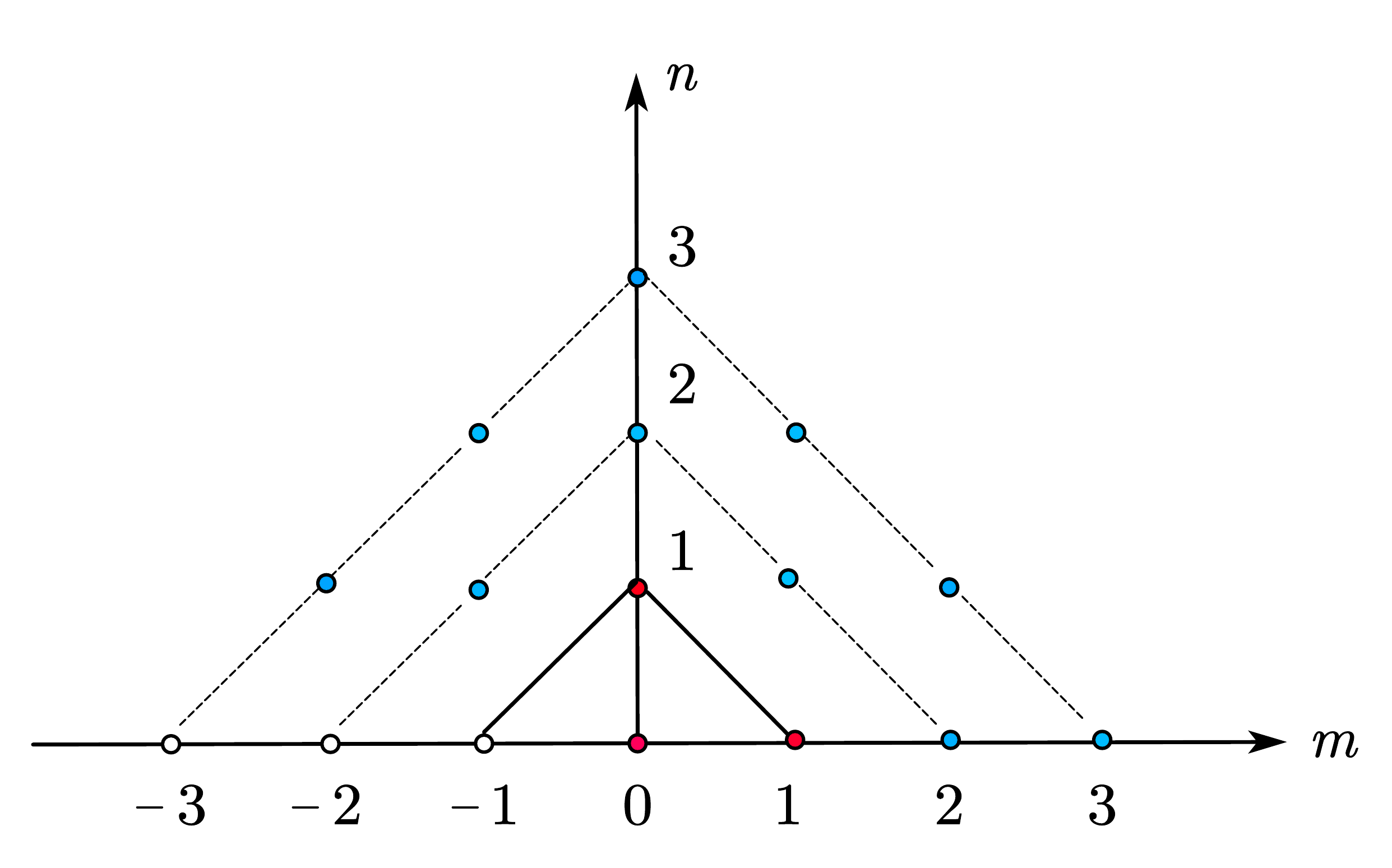}
	 \caption{Illustration of the 2D Fourier series and its cubic expansion.}
	 \label{FIG:2DFourier}
\end{figure}

Before the proof, we would introduce a property of the 2D Fourier series. Figure \ref{FIG:2DFourier} shows the relation between 2D Fourier truncation ($p=1$) and its high-order harmonic terms. Notice that a filled circle represents a harmonic term included in computing, and a hollow one does not exist. Truncation harmonics are geometrically denoted as coordinates by the parameter pairs $(m,n)$, and all these coordinates constitute a point set $\Omega _{p}$ (see red filled points in Figure \ref{FIG:2DFourier}). As provided in Appendix A, high-order terms not in Eq. (\ref{eq:represents}) are generated by cubic expansion, and their parameters constitute another set $\Omega '_p$ (see blue points in Figure \ref{FIG:2DFourier}). We conclude that the polynomial nonlinear function $x^{\phi}$ ($\phi$ refers to the nonlinearity) actually enlarges the original point set $\Omega _{p}$ to a bigger one $\Omega _{\phi p}=\Omega _p\cup \Omega '_p$.

\begin{theorem}[Frequency ratio is irrational]\label{thm1}
Suppose a system with nonlinearity $\phi$, the response has two base frequencies and the ratio  $\omega_1/\omega_2$ is irrational. The RMHB is equivalent to the MHB only if sampling period $T$ and the number of collocations $M$ are both infinite.
\end{theorem}

\begin{proof}
Comparing algebraic Eqs. (\ref{MHB NAEs}) and (\ref{RMHB NAEs}) derived by two methods respectively, they are equivalent when meeting (a), $\mathbf{E}^*\mathbf{E}=\mathbf{I}_{N}\otimes \mathbf{A}$ is an identity matrix, and (b), $\mathbf{\hat{f}}(\mathbf{\hat{x}})=\mathbf{E}^*\mathbf{\tilde{f}}_M(\mathbf{E\hat{x}})$.

\begin{enumerate}
    \item Each element in matrix $\mathbf{A}$ can be written in trigonometric summations. It will convert into the integral, as $M$ increases to infinity. For elements that row and column index $i\ne j$ in the  matrix $\mathbf{A}$, e.g.,
\begin{flalign}
    \frac{2}{M} \lim _{M \rightarrow \infty} \sum_{i}^{M} \sin \left(a \omega_{1}+b \omega_{2}\right) t_{i} \cdot \cos \left(c \omega_{1}+d \omega_{2}\right) t_{i} \\
    =\frac{2}{T} \int_{0}^{T} \sin \left(a \omega_{1}+b \omega_{2}\right) t \cdot \cos \left(c \omega_{1}+d \omega_{2}\right) t \mathrm{~d} t, \nonumber
\end{flalign}

parameter pairs $(a,b),(c,d)$ satisfy the inequality relation (\ref{relation}). Let $\theta _1=a\omega _1+b\omega _2$, $\theta _2=c\omega _1+d\omega _2$ and apply variable substitution $\tau =\theta _2t/T$. Since function $f\left( \tau \right) =\sin (\frac{\theta _1}{\theta _2}T\tau )$ is Riemann integrable over the integral interval, applying the Riemann-Lebesgue lemma \cite{Serov2017}:
\begin{align*}
    \lim _{T \rightarrow\infty} \frac{2}{T} \int_{0}^{T} \sin \theta_{1} t \cdot \cos \theta_{2} t \mathrm{~d} t \stackrel{\tau=\frac{\theta_{2} t}{T}} {=} &\lim _{T \rightarrow \infty} \frac{2}{\theta_{2}} \int_{0}^{\theta_{2}} \sin \left(\frac{\theta_{1} T}{\theta_{2}} \tau\right) \cos (T \tau) \mathrm{d} \tau \\
    =&\lim_{T\rightarrow \infty} \frac{2}{\theta _2}\int_0^{\theta _2}{f\left( \tau \right) \cdot \cos \left( T\tau \right)}\mathrm{d}\tau =0,\nonumber
\end{align*}

the value of the above integral is 0, other off-diagonal elements can be proved similarly. For those elements on the diagonal, e.g.,
\begin{equation}
    \frac{2}{M}\lim_{M\rightarrow \infty} \sum\nolimits_i^M{\cos \left( \theta t_i \right)}\cdot \cos \left( \theta t_i \right) =\frac{2}{T}\int_0^T{\left( \cos \theta t \right) ^2\mathrm{d}t},
    \label{eq:diagonal}
\end{equation}
since the integrand function $g=\cos ^2\theta $ does not change sign on $[0,T]$, suppose $f(\tau)=1$ is the continuous function, exploit an integral limit conclusion
\begin{equation}
    \lim_{n\rightarrow \infty} \int_a^b{f\left( x \right) g\left( nx \right) \mathrm{d}x=\left( \int_a^b{f\left( x \right) \mathrm{d}x} \right)}\left( \frac{1}{T}\int_0^T{g\left( x \right) \mathrm{d}x} \right),\nonumber
\end{equation}
the Eq. (\ref{eq:diagonal}) is 1, we can also prove other elements on the diagonal in $\mathbf{A}$.
\item For brevity, here we assume a single DOF system, $\mathbf{\hat{f}}(\mathbf{\hat{x}})$ and $\mathbf{\tilde{f}}_M(\mathbf{E\hat{x}})$ are assumed as
$\mathbf{\hat{f}}(\mathbf{\hat{x}})=[\hat{f}_{c}(0,0) \,\, \hat{f}_{c}(1,0) \,\, \hat{f}_{s}(1,0)\,\, \cdots \,\, \hat{f}_{c}(0,p)\,\, \hat{f}_{s}(0,p)]^{\mathrm{T}}$ and $\mathbf{\tilde{f}}_M\left( \mathbf{E\hat{x}} \right) =[\mathbf{f}\left( \mathbf{x}\left( t_1 \right) \right) \,\, \mathbf{f}\left( \mathbf{x}\left( t_2 \right) \right)\,\, \cdots \\ \,\, \mathbf{f}\left( \mathbf{x}(t_M) \right) ]^{\mathrm{T}}$, where symbols with hat like $\hat{f}_{c}(1,0)$ denote the Fourier coefficients, and $\mathbf{f}(\mathbf{x}(t_i))$ is the corresponding temporal quantity at prescribed time instant $t_i$,
\begin{equation}
    \mathbf{f}\left( \mathbf{x}\left( t_i \right) \right) =\sum_{\left( \bar{m},\bar{n} \right) \in \Omega _{\phi p}}[\hat{f}_c\left( \bar{m},\bar{n} \right) \mathrm{c}^{\bar{m},\bar{n}}(t_i)+\hat{f}_s\left( \bar{m},\bar{n} \right) \mathrm{s}^{\bar{m},\bar{n}}(t_i)].\nonumber
\end{equation}

For nonlinearity $\phi=1$, the $\mathbf{\hat{f}}(\mathbf{\hat{x}})=\mathbf{E}^*\mathbf{\tilde{f}}_M(\mathbf{E\hat{x}})$ holds.

For $\phi \geqslant 2$, $\mathbf{\tilde{f}}_M(\mathbf{E\hat{x}})$ can be divided into two parts: $p$-order truncation harmonics and high-order terms as

\begin{equation}
\begin{aligned}
    \tilde{\mathbf{f}}_{M}=\left[\begin{array}{c}
    \sum_{(m, n) \in \Omega_{p}} \hat{f}_{c}(m, n) \mathrm{c}^{m, n}\left(t_{1}\right)+\hat{f}_{s}(m, n) \mathrm{s}^{m, n}\left(t_{1}\right) \\
    \sum_{(m, n) \in \Omega_{p}} \hat{f}_{c}(m, n) \mathrm{c}^{m,     n}\left(t_{2}\right)+\hat{f}_{s}(m, n) \mathrm{s}^{m, n}\left(t_{2}\right) \\
\vdots \\
    \sum_{(m, n) \in \Omega_{p}} \hat{f}_{c}(m, n) \mathrm{c}^{m, n}\left(t_{M}\right)+\hat{f}_{s}(m, n) \mathrm{s}^{m, n}\left(t_{M}\right)
    \end{array}\right]
\end{aligned}\nonumber
\end{equation}

\begin{equation}
\quad +\left[ \begin{array}{c}
	\sum_{\left( m',n' \right) \in \varOmega _{p}^{'}}{\hat{f}_c\left( m',n' \right) \mathrm{c}^{m',n'}\left( t_1 \right) +\hat{f}_s\left( m',n' \right) \mathrm{s}^{m',n'}\left( t_1 \right)}\\
	\sum_{\left( m',n' \right) \in \varOmega _{p}^{'}}{\hat{f}_c\left( m',n' \right) \mathrm{c}^{m',n'}\left( t_2 \right) +\hat{f}_s\left( m',n' \right) \mathrm{s}^{m',n'}\left( t_2 \right)}\\
	\vdots\\
	\sum_{\left( m',n' \right) \in \varOmega _{p}^{'}}{\hat{f}_c\left( m',n' \right) \mathrm{c}^{m',n'}\left( t_M \right) +\hat{f}_s\left( m',n' \right) \mathrm{s}^{m',n'}\left( t_M \right)}\\
\end{array} \right] . \nonumber
\end{equation}

The first term of the above formula is expressed as $\mathbf{E}\mathbf{\hat{f}}(\mathbf{\hat{x}})$. Because $\mathbf{E}^* \mathbf{E}$ is an identity matrix, multiply transformation matrix $\mathbf{E}^*$ on both sides of the above formula, we derive
\begin{equation}
    \mathbf{E}^*\mathbf{\tilde{f}}_{M}\left( \mathbf{E\hat{x}} \right) =\mathbf{E}^*\mathbf{E\hat{f}}\left( \mathbf{\hat{x}} \right) +\mathbf{E}^*\mathbf{E}_1\mathbf{\hat{f}}'\left( \mathbf{\hat{x}} \right) =\mathbf{\hat{f}}\left( \mathbf{\hat{x}} \right) +\mathbf{E}^*\mathbf{E}_1\mathbf{\hat{f}}'\left( \mathbf{\hat{x}} \right) ,\label{eq:aliasing}
\end{equation}
denote
\begin{equation}
    \mathbf{E}_{1}=\left[\begin{array}{c}
    \mathbf{F}_{h}\left(t_{1}\right) \\
    \mathbf{F}_{h}\left(t_{2}\right) \\
    \vdots \\
    \mathbf{F}_{h}\left(t_{M}\right)
    \end{array}\right], \quad 
    \hat{\mathbf{f}}^{\prime}(\hat{\mathbf{x}})=\left[\begin{array}{c}
    \hat{f}_{c}(p+1,0) \\
    \hat{f}_{s}(p+1,0) \\
    \vdots \\
    \hat{f}_{c}\left(m^{\prime}, n^{\prime}\right) \\
    \hat{f}_{s}\left(m^{\prime}, n^{\prime}\right) \\
    \vdots \\
    \hat{f}_{c}(0, \phi p) \\
    \hat{f}_{s}(0, \phi p)
    \end{array}\right],\quad
    \mathbf{F}_{h}\left(t_{i}\right)=\left[\begin{array}{c}
    \cos (p+1) \omega_{1} t_{i} \\
    \sin (p+1) \omega_{1} t_{i} \\
    \vdots \\
    \cos \left(m^{\prime} \omega_{1}+n^{\prime} \omega_{2}\right) t_{i} \\
    \sin \left(m^{\prime} \omega_{1}+n^{\prime} \omega_{2}\right) t_{i} \\
    \vdots \\
    \cos \phi p \omega_{2} t_{i} \\
    \sin \phi p \omega_{2} t_{i}
    \end{array}\right]^{\mathrm{T}}.\nonumber
\end{equation}
Where $ i=1,2, \ldots, M$. Define $\mathbf{E}_{\mathrm{A}}=\mathbf{E}^*\mathbf{E}_1$ as "aliasing matrix" \cite{Dai2022}. Each element in the matrix $\mathbf{E}_{\mathrm{A}}$ can be written as the integral limit, with $M$ and $T$ being both infinite. By applying the Riemann-Lebesgue lemma, it can be proved that all elements in the $\mathbf{E}_{\mathrm{A}}$ gradually go to zero. Only satisfying the sampling rule promises conditional equivalence.
\end{enumerate}

\end{proof}

One completely different kind of multiple harmonic balance problem is periodic but with some relatively high-frequency components. Another sampling rule will be given based on the proposed RMHB method.
\begin{theorem}[Frequency ratio is rational]\label{thm2}
Suppose a system with nonlinearity $\phi$, the response has two base frequencies and the ratio  $\omega_1/\omega_2$ is rational. The RMHB is equivalent to the MHB only if sampling period $T=2{{\pi}/{\mathrm{G}}}\mathrm{CD(}\omega _1,\omega _2)$ and the number of collocations
\begin{eqnarray}
    M>(\phi +1)\frac{p\cdot \max \left( \omega _1,\omega _2 \right)}{\mathrm{GCD}\left( \omega _1,\omega _2 \right)},
\end{eqnarray}
where $\mathrm{GCD}$ is the greatest common divisor of two numbers.
\end{theorem}

\begin{proof}
    The sampling period $T$ can be determined by the common period of the two input frequencies \cite{Stupel2012}. It is similar to the proof of Theorem \ref{thm1}, the RMHB and the MHB are equivalent when both conditions hold: (a), $\mathbf{E}^*\mathbf{E}=\mathbf{I}_{N}\otimes \mathbf{A}$, and (b), $\mathbf{\hat{f}}(\mathbf{\hat{x}})=\mathbf{E}^*\mathbf{\tilde{f}}_M(\mathbf{E\hat{x}})$. In order to assist our derivation process, note the following well-known property of trigonometric functions \cite{Bocher1906,Jackson1913} that if $x_1,\cdots,x_M$ are $M$ nodes disposed of successive time interval of $2\pi/M$, in which $\alpha,\beta,\alpha+\beta,\alpha-\beta$ are positive integers but less than number of collocations $M$.
\begin{equation}
    \left\{\begin{array}{l}
    \sum_{i=1}^{i=M} \sin \left(\alpha x_{i}\right) \cos \left(\beta x_{i}\right)=0, \\
    \sum_{i=1}^{i=M} \cos \left(\alpha x_{i}\right) \cos \left(\beta x_{i}\right)=\begin{cases}0 & (\alpha \neq \beta), \\
    \frac{M}{2} & (\alpha=\beta),\end{cases} \\
    \sum_{i=1}^{i=M} \sin \left(\alpha x_{i}\right) \sin \left(\beta x_{i}\right)= \begin{cases}0 & (\alpha \neq \beta), \\
    \frac{M}{2} & (\alpha=\beta) .\end{cases}
\end{array}\right.\label{eq:trig}
\end{equation}

\begin{enumerate}
    \item Let $t_i=\frac{T}{M}i, x_i=\frac{2\pi}{M}i$,  each element in matrix $\mathbf{A}$ can be written as a trigonometric summation, e.g.,
    \begin{equation}
        \sum_i^M{\sin}\left( a\omega _1+b\omega _2 \right) t_i\cdot \cos \left( c\omega _1+d\omega _2 \right) t_i=\sum_i^M{\sin}\left( a\omega _1+b\omega _2 \right) \frac{T}{2\pi}x_i\cdot \cos \left( c\omega_1+d\omega _2 \right) \frac{T}{2\pi}x_i,\nonumber
    \end{equation}
    with $(a,b),(c,d)\in \Omega _p$. Period $T=2{{\pi}/{\mathrm{G}}}\mathrm{CD(}\omega _1,\omega _2)$ can make that $\left( a\omega _1+b\omega _2 \right) \frac{T}{2\pi}$, $\left( c\omega _1+d\omega _2 \right) \frac{T}{2\pi}$ and $\left( \left( a\pm c \right) \omega _1+\left( b\pm d \right) \omega _2 \right) \frac{T}{2\pi}$ are all integers. Applying (\ref{eq:trig}) for each element in $\mathbf{A}$,  $\mathbf{E}^{*}\mathbf{E}$ is an identity matrix when the number of collocations $M$ 

    \begin{align*}
        M>(\phi +1)\frac{p\cdot \max \left( \omega _1,\omega _2 \right)}{\mathrm{GCD}\left( \omega _1,\omega _2 \right)}& \geqslant \frac{2p\cdot \max \left( \omega _1,\omega _2 \right)}{\mathrm{GCD}\left( \omega _1,\omega _2 \right)}\\
        & \geqslant \max \left[ a\omega _1+b\omega _2,c\omega _1+d\omega _2,\left( a\pm c \right) \omega _1+\left( b\pm d \right) \omega _2 \right] \frac{T}{2\pi}.\nonumber
    \end{align*}
        
    \item According to the definition of the 2D Fourier series, the maximum harmonic contained in the $p$-order RMHB is $p\cdot \max(\omega_{1},\omega_{2})$, and the maximum component obtained by the nonlinear function $x^{\phi}$ is $\phi p\cdot \max(\omega_{1},\omega_{2})$. Here all elements in the aliasing matrix $\mathbf{E}_{\mathrm{A}}$ can be written in the form of a trigonometric summation, it can be a zero matrix, when
    \begin{equation}
        \begin{aligned}
            M>&(\phi +1)\frac{p\cdot \max \left( \omega _1,\omega _2 \right)}{\mathrm{GCD}\left( \omega _1,\omega _2 \right)}=\left( \left( \phi +1 \right) p\cdot \max \left( \omega _1,\omega _2 \right) \right) \frac{T}{2\pi}\\
           &\geqslant \max \left( m\omega _1+n\omega _2,m'\omega _1+n'\omega _2,\left( m\pm m' \right) \omega _1+\left( n\pm n' \right) \omega _2 \right) \frac{T}{2\pi},\nonumber
        \end{aligned}
    \end{equation}
    
\end{enumerate}
\end{proof}

\begin{remark}
Suppose a system with nonlinearity $\phi$, the response has multiple base frequencies $\left\{ \omega _1,\omega _2,\cdots ,\omega _t \right\} $ and all ratios are rational. The RMHB is equivalent to the MHB only if sampling period $T=2{{\pi}/{\mathrm{G}}}\mathrm{CD(}\omega _1,\omega _2,\cdots ,\omega _t)$ and the number of collocation nodes
\begin{equation}
    M>(\phi +1)\frac{p\cdot \max \left( \omega _1,\omega _2,...,\omega _t \right)}{\mathrm{GCD}\left( \omega _1,\omega _2,...,\omega _t \right)}.
\end{equation}
\end{remark}

To research the aliasing phenomenon in the multiple harmonic balance computations, we consider a Duffing equation with two input signals \cite{Prabith2020}
\begin{equation}
    \ddot{x}+0.2\dot{x}+x+0.2x^3=3\cos 4t+5\cos 2.8t.\nonumber
\end{equation}

The critical number of collocation nodes $M=41$ for de-aliasing can be given by Theorem \ref{thm2}. To intuitively demonstrate the effect of adding nodes to eliminate aliasing, the Monte-Carlo simulation results are represented. This research applies the RMHB1 and all frequency unknowns are specified within the range of $\pm 5$ unites, we select 10000 groups of such random initials for simulation. For simplicity, the $p$-th order RMHB method, i.e., the RMHB with $p$-order multidimensional Fourier truncation, is denoted as RMHB$p$.

\begin{figure}
	\centering
	 \includegraphics[scale=0.4]{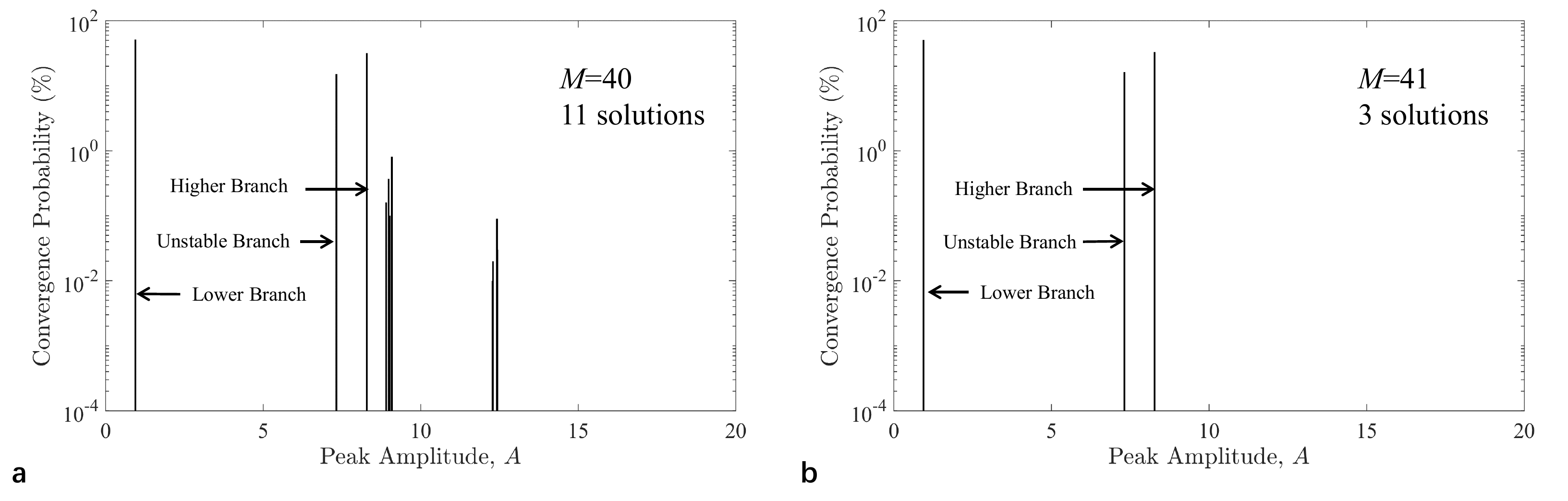}
	 \caption{Comparison of Monte-Carlo histogram for the probability distribution of solutions for $\ddot{x}+0.2\dot{x}+x+0.2x^3=3\cos 4t+5\cos 2.8t$ with different number of collocations (a), $M=40$ and (b), $M=41$.}
	 \label{FIG:MonteCarlo}
\end{figure}

Figure \ref{FIG:MonteCarlo} indicates there are three solutions to the system: the higher branch and the lower branch are stable,  and one unstable branch (cannot be captured by numerical integration). Table \ref{tabM} shows that the RMHB method obtains fewer non-physical solutions as $M$ increases. When $M=20$, probabilities corresponding to the upper, lower, and unstable branches are $46.52\%$, $10.79\%$, and $18.42\%$ respectively,  while $24.27\%$ is the probability of non-physical solutions. As gradually add the $M$ to the critical number, the RMHB method will only get 3 physical solutions with proper sampling, which accounts for $53.45\%$, $15.01\%$, and $31.54\%$ probability when the number of collocations meets the condition. 

\begin{figure}
	\centering
	 \includegraphics[scale=0.4]{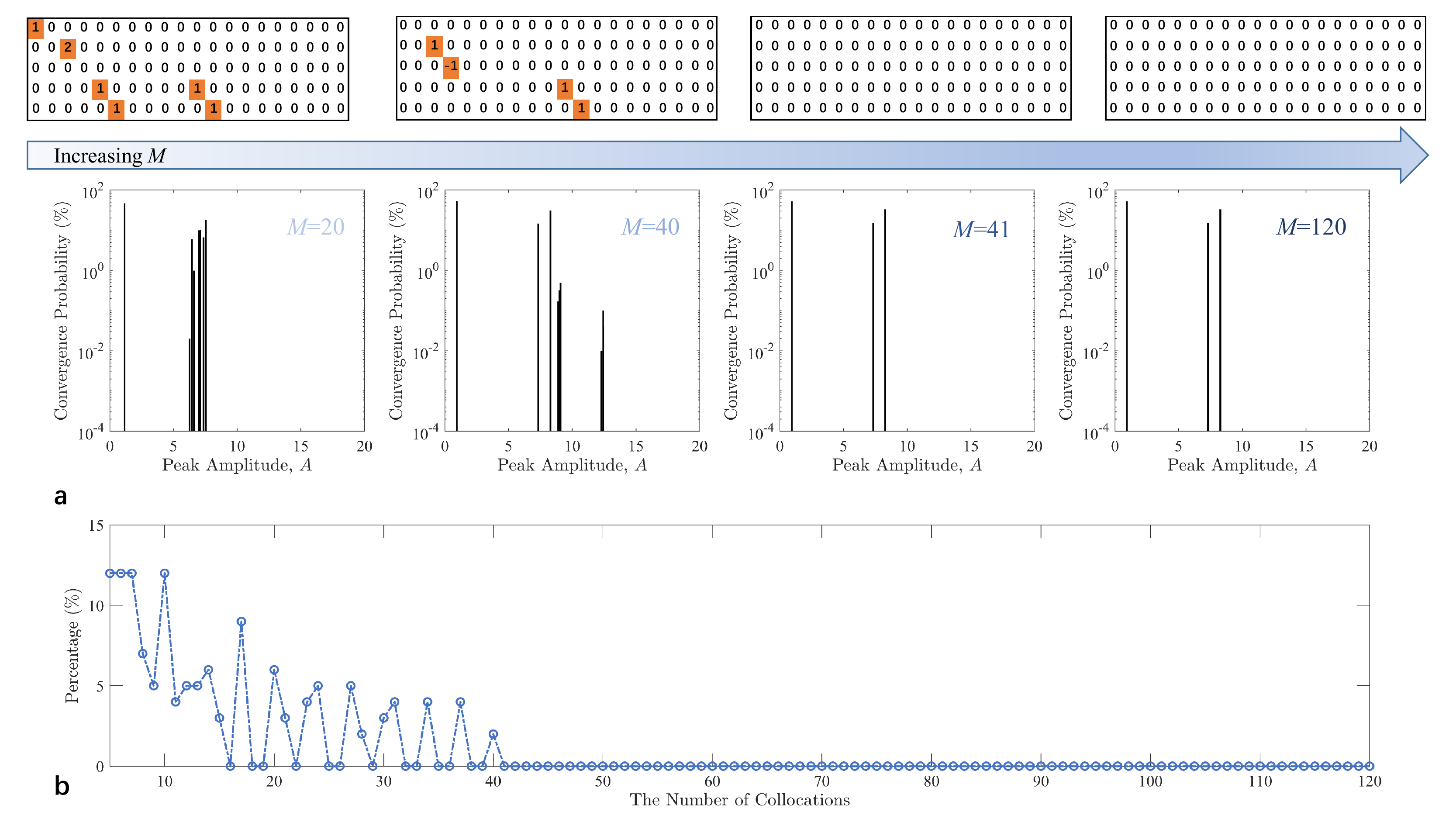}
	 \caption{Illustration of the aliasing matrix changing with the number of collocations $M$ from 5 to 120 when solving $\ddot{x}+0.2\dot{x}+x+0.2x^3=3\cos 4t+5\cos 2.8t$. (a), Change of the aliasing matrix and its effect to eliminate nonphysical solutions. (b), Change of percentage of non-zero elements in the aliasing matrix.}
	 \label{FIG:Aliasing matrix}
\end{figure}

Because the aliasing can be quantified as whether the matrix $\mathbf{E}_{\mathrm{A}}$ is a zero matrix or not. Figure \ref{FIG:Aliasing matrix}a demonstrates the dynamic change of the aliasing matrix state and the de-aliasing effect. Besides, Figue \ref{FIG:Aliasing matrix}b indicates that Theorem \ref{thm2} is a sufficient condition for eliminating aliasing, i.e., not all the $M$ less than the critical value will cause aliasing.

${p\cdot \max \left( \omega _1,\omega _2 \right)}/{\mathrm{GCD}\left( \omega _1,\omega _2 \right)}$ is comparable to the truncation order $N$ in the RHB method (the base frequency is $\mathrm{GCD(}\omega _1,\omega _2)$), thus Theorem \ref{thm2} can be rewritten in a more familiar form: $M>(\phi+1)N$ \cite{Dai2022}. Herein we tell that the RMHB method compresses redundant frequency variables by introducing multiple base frequencies, but the number of collocations $M$ required for de-aliasing is consistent with the single-frequency method. So the percentage of non-zero elements in the matrix $\mathbf{E}_{\mathrm{A}}$ varies intermittently. This phenomenon only exists in the case of multiple base frequencies and is not regular, depending on the actual response and the truncation. However, the critical value in Theorem \ref{thm2} is uniquely determined, and the equivalence holds only when $M$ exceeds the critical value.

\begin{table}[htbp]
\centering
\footnotesize{
 \caption{\label{tabM}Solution details for various collocation settings}
 \begin{tabular}{ccccc}
  \toprule
    $M$ & Lower Branch (\%) & Unstable Branch (\%) & Higher Branch (\%) & Non-physical (\%) \\ %
  \midrule
 20 & 46.52 & 10.79 & 18.42 & 24.27 \\
 40 & 53.10 & 14.20 & 31.30 & 1.40 \\
 41 & 53.45 & 15.01 & 31.54 & 0 \\
 60 & 53.38 & 15.01 & 31.61 & 0 \\
  \bottomrule
 \end{tabular}}
\end{table}

\section{Results}
\subsection{Forced Self-Excited Van der Pol Equation}
The Van der Pol equation is 
\begin{equation}
    \ddot{x}-\varepsilon \left( 1-x^2 \right) \dot{x}+x=F\cos \omega _1t,\label{eq:VanderPol}
\end{equation}
where damping coefficient $\varepsilon=0.1$, amplitude $F=0.25$. The forcing frequency $\omega_{1}=4/\pi$ is near the primary resonance 1.0, the ratio of forcing frequency $\omega_{1}$ and the primary resonance $\omega_{2}=1$ is irrational. Especially such a specific system coexists with two different base frequencies and their linear combinations, the response is not periodic and this kind of motion is also described as "mild" chaos \cite{Liu2007}. By recasting \cite{Cochelin2009} 
\begin{equation}
    \begin{cases}
	\dot{x}=u,\\
	\dot{u}=\left( \varepsilon u-x \right) -\varepsilon x^2u+F\cos \omega _1t,\\
    \end{cases} \nonumber
\end{equation}
the original problem  (\ref{eq:VanderPol}) transforms to a typical polynomial system with $\phi=3$. The algebraic equations for determining harmonic components can be written as
\begin{equation}
    \begin{cases}
	\mathbf{E}^*\mathbf{E} \nabla \mathbf{\hat{x}}=\mathbf{E}^*\mathbf{E\hat{u}},\\
	\mathbf{E}^*\mathbf{E} \nabla \mathbf{\hat{u}}=\mathbf{E}^*\left( \varepsilon \mathbf{E\hat{u}}-\mathbf{E\hat{x}} \right) -\varepsilon \mathbf{E}^*\tilde{\mathbf{f}}_M +\mathbf{E}^*\mathbf{\tilde{H}},\\
\end{cases}
\end{equation}
with $\tilde{\mathbf{f}}_M=\left[\begin{array}{c}
    x^2\left(t_1\right) u\left(t_1\right) \\
    x^2\left(t_2\right) u\left(t_2\right) \\
    \vdots \\
    x^2\left(t_M\right) u\left(t_M\right)
    \end{array}\right]$, \,\,
    $\tilde{\mathbf{H}}=\left[\begin{array}{c}
    F \cos \left(\omega_1 t_1\right) \\
    F \cos \left(\omega_1 t_2\right) \\
    \vdots \\
    F \cos \left(\omega_1 t_M\right)
    \end{array}\right]$.

When using the HDHB method with two base frequencies to compute the quasi-periodic response, Liu \cite{Liu2007} firstly pointed out that proper sampling period $T$ and the number of collocations $M$ are equally crucial for calculation accuracy. This discovery also applies to the RMHB method. Table \ref{tab1} shows that only when $T$ and $M$ are sufficient, the aliasing error can be controlled. Then the RMHB method will produce almost the same result as the MHB method.

\begin{figure}
	\centering
	 \includegraphics[scale=0.4]{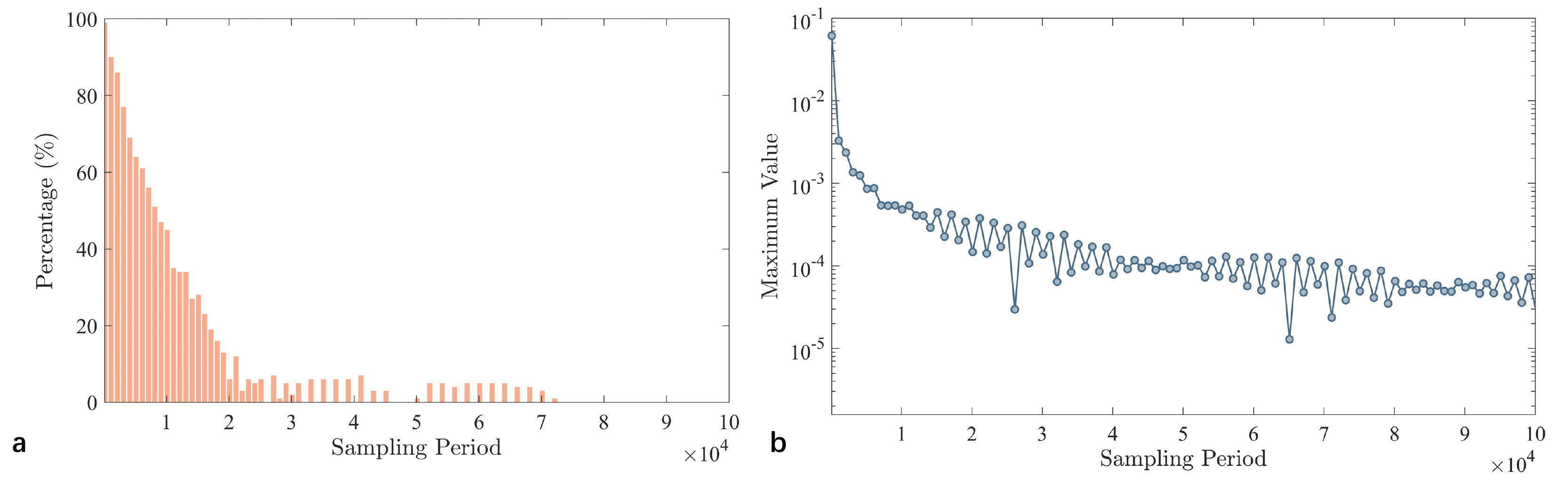}
	 \caption{Effect of $T$ and $M$ on aliasing matrix properties. (a), Percentage of elements in aliasing matrix with absolute value greater than $10^{-4}$. (b), The change of maximum element value in the aliasing matrix}
	 \label{FIG:Ration & Max}
\end{figure}

Next, we start from a statistical view to demonstrate the improvement effect on the properties of the aliasing matrix $\mathbf{E}_{A}$ with enough $T$ and $M$. Set each collocation node in the unit time,  Figure \ref{FIG:Ration & Max} shows that increasing $T$ and $M$ can reduce both the magnitude of the maximum element and the percentage of elements greater than $10^{-4}$. When $T$ and $M$ reach $7\times10^4$, all elements in the aliasing matrix will no longer exceed $10^{-4}$. From Table \ref{tab1} and Figure \ref{FIG:Ration & Max}, we conclude that sufficient $M$ and $T$ are necessary to eliminate aliasing.

\begin{figure}
	\centering
	 \includegraphics[scale=0.25]{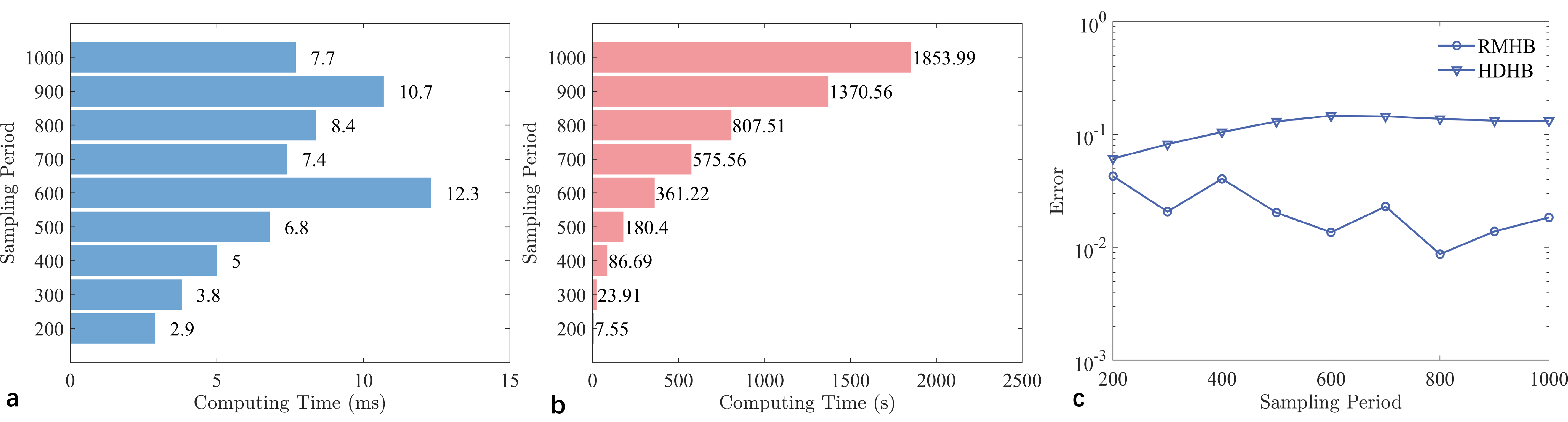}
	 \caption{Comparison of computational efficiency and accuracy between the RMHB and HDHB methods. (a, b), Computing time of the RMHB method (a) and HDHB method (b) for solving forced Van der Pol oscillator. (c), The computational error of the two methods.}
	 \label{FIG:RMHB & HDHB}
\end{figure}

Since the number of collocations is variable, the transformation matrix $\mathbf{E}^*$ was regarded as the pseudo inverse of the collocation matrix $\mathbf{E}$ \cite{Chua1981,Liu2007,Lindblad2022}. Thus the HDHB method is is similar in nature to the time domain collocation (TDC) method \cite{Dai2012}. Figure \ref{FIG:RMHB & HDHB} shows the total computing time and error curve for the RMHB1 and the HDHB method (same harmonic components) by using several sets of sampling periods (interpolation nodes set per unit time). The reference solution is obtained by the RK4 method. We found that the computing time of the RMHB method with the same order is less affected by the $T$ and $M$, alleviates the computational burden and the whole time does not exceed 15 milliseconds at most. However the HDHB method is a time domain approach, variables in the NAEs are numerical values on the time domain nodes \cite{Liu2006,Labryer2009}, adding $M$ will inevitably heavier computational burden. Figure \ref{FIG:RMHB & HDHB}b reveals that the computing time of the HDHB method reaches nearly 30 minutes when $T=M=1000$. Besides Figure \ref{FIG:RMHB & HDHB}c shows that the error of the RMHB1 tends to converge as the $T$ and $M$ are sufficient. But the error of the HDHB method remains in the order of $10^{-1}$, because the transformation matrix $\mathbf{E}^{*}$ is derived by numerical calculation and no longer maintains physical meaning. Therefore $T$ and $M$ have a limited contribution to improving the convergence in the HDHB computation. In short, compared with the RMHB method, the computational efficiency of the HDHB method is limited by the number of collocation nodes, and the accuracy cannot be guaranteed.

\begin{table}[htbp]
\centering
\footnotesize{
 \caption{\label{tab1}Results of the RMHB1 for different $T$ and $M$}
 \begin{tabular}{ccccc}
  \toprule
    $T$  & $M$ & $\sqrt{x_{1,0}^{2}+y_{1,0}^{2}}$ & $\sqrt{x_{0,1}^{2}+y_{0,1}^{2}}$ & Amplitude error \\ 
  \midrule
    100 & 100 & 0.4112 & 2.004 & 0.1089 \\
    500 & 500 & 0.3952 & 1.934 & 0.0175\\
    5000 & 5000 & 0.3960 & 1.922 & 0.0022 \\
    $5\times 10^{4}$ & $5\times 10^{4}$ & 0.3961 & 1.920 & $2.4\times10^{-4}$ \\
    $5\times 10^{5}$ & $5\times 10^{5}$ & 0.3960 & 1.920 & $2.6\times10^{-5}$ \\
  \bottomrule
 \end{tabular}}
\end{table}

\begin{figure}
	\centering
	 \includegraphics[scale=0.35]{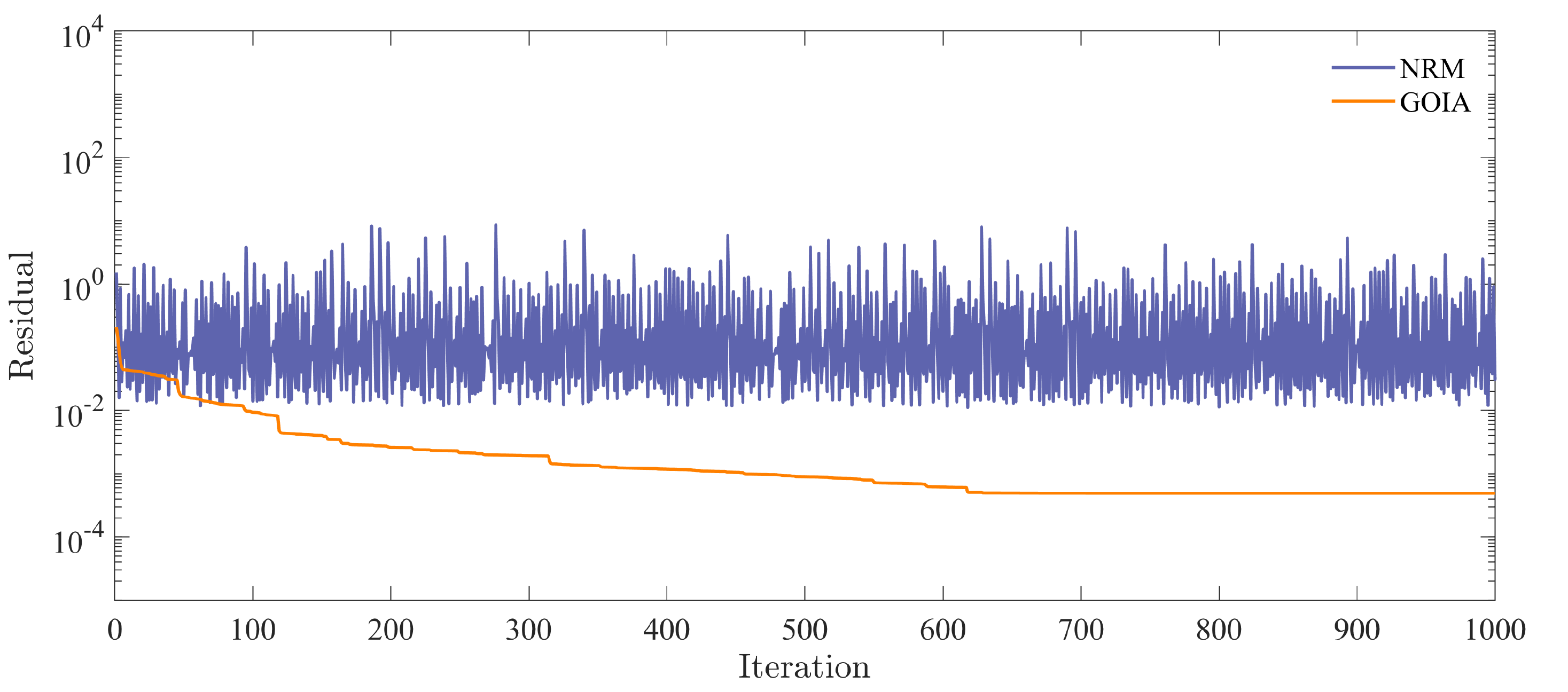}
	 \caption{Residual error versus iterative curve for solving the RMHB3 algebraic equations by using two kinds of NAEs solvers}
	 \label{FIG:NRM & GOIA}
\end{figure}

The RMHB method can greatly reduce the symbolic operations, but the higher-order harmonic estimation may suffer from the ill-conditioned problem when solving quasi-periodic response. In the practical numerical calculation, only a relatively large finite value can be selected as $T$ and $M$, resulting in differences in some coefficients of the NAEs, which brings difficulties to the Newton method. Taking RMHB3 as an example, our research utilizes the Newton-Raphson method (NRM) and global optimal iterative algorithm (GOIA) \cite{Liu2012} respectively. It is worth mentioning that, different from the NRM, GOIA is one kind of scalar homotopy method. GOIA works by finding the best descend vector to iteratively solve a system of NAEs $\mathbf{F}(\mathbf{x})=0$, without requiring the inversion of the Jacobian matrix. The iteration form can be written in

\begin{equation}
    \mathbf{x}_{k+1}=\mathbf{x}_k-\left( 1-\gamma \right) \frac{\mathbf{F}_k\cdot \mathbf{v}_k}{\left\| \mathbf{v}_k \right\| ^2}\mathbf{u}_k,
\end{equation}
where $0\leqslant \gamma <1$, $\mathbf{F}_{k}=\mathbf{F}(\mathbf{x}_k)$ is the iteration residual for each step, detailed calculation of the optimal descent vector $\mathbf{u}_k$ and combination vector $\mathbf{v}_k$ can refer to \cite{Liu2012}. Set $T=M=2\times10^4$ and initial value are $\hat{x}_{c}(1,0)=0.3$,  $\hat{x}_{c}(0,1)=1.9$. Figure \ref{FIG:NRM & GOIA} shows that it is necessary to select an appropriate NAEs solver, the NRM is fast and efficient, but depends on the initial values, and may not converge because of ill-conditioned problems. Some optimal methods like the GOIA method and the Tikhonov regularization method \cite{Wang2020,Zheng2021} can overcome the above-mentioned iterative convergence problem but it is time-consuming. 

\begin{figure}
	\centering
	 \includegraphics[scale=0.25]{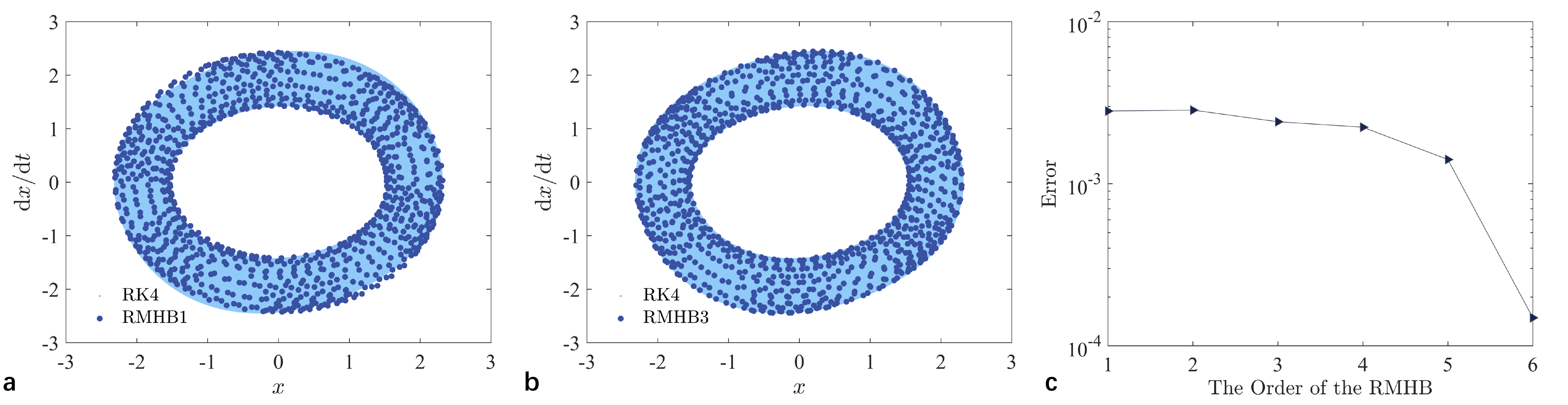}
	 \caption{(a, b), Phase plot calculated by RMHB method with $p=1$ (a) and $p=3$ (b). (c), Computing errors vary with the adopted order of the RMHB method.}
	 \label{FIG:p1 & p3}
\end{figure}

Compared with the numerical integration, the phase plot and the error curve obtained by the RMHB method with different orders are shown in Figure \ref{FIG:p1 & p3}. Set $T=M=10^5$, we can clearly tell the convergence tendency of the RMHB in Figure \ref{FIG:p1 & p3}c when solving quasi-periodic responses. More specifically, from Table \ref{tabc} we conclude that sufficient sampling periods, collocations, and harmonic components can contribute to the computing accuracy. Let truncation order $p=5$, $M=T=10^5$, the RMHB method controls the amplitude error to $10^{-4}$ when solving quasi-periodic response (base frequencies are incommensurable).

\begin{table}[htbp]
\centering
\footnotesize{
 \caption{\label{tabc}Comparison of amplitude error of the RMHB method with different orders, sampling periods, and collocations}
 \begin{tabular}{cccc}
  \toprule
    Order & $M$ & $T$ & Amplitude error\\ 
  \midrule
     1 & $10^3$ & $10^3$ & 0.0177 \\
    1 & $10^5$ & $10^5$ & 0.0028 \\
    3 & $10^3$ & $10^3$ & 0.0058 \\
    3 & $10^5$ & $10^5$ & 0.0024 \\
    5 & $10^5$ & $10^5$ & 0.0014 \\
    6 & $10^5$ & $10^5$ & $1.49 \times 10^{-4}$ \\
  \bottomrule
 \end{tabular}}
\end{table}

\subsection{Duffing Equation with Two-Frequency Inputs}
The state equation of Duffing oscillator  is
\begin{equation}
    \begin{cases}
	\dot{x}_1=x_2,\\
	\dot{x}_2=-cx_2-kx_1-\alpha x_{1}^{3}+y\left( t \right),\\
\end{cases}\label{eq:duffing}
\end{equation}
and $y(t)$ is a two-frequency input signal
\begin{equation}
    y(t)=A_1\cos \omega_1+A_2\cos \omega_2. 
\end{equation}
Eq. (\ref{eq:duffing}) can commonly describe dynamics models (eg., ferroresonance circuits, differential-pair amplitude modulator circuits) and show unique physical phenomena like sub-harmonic, quasi-periodic, and chaos solutions.

\begin{figure}
	\centering
	 \includegraphics[scale=0.4]{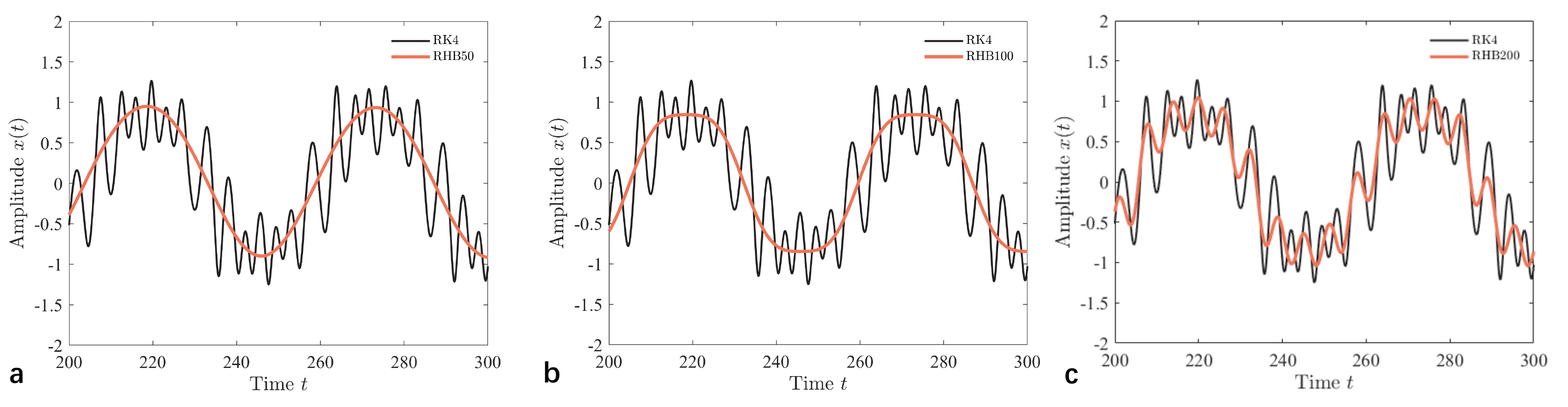}
	 \caption{Comparison of the single base frequency RHB method for solving $\ddot{x}+0.05\dot{x}+x+x^3=0.3\cos t+1.5\cos 0.115t$ with (a), $N=50$, (b), $N=100$, and (c), $N=200$.}
	 \label{FIG:RHB}
\end{figure}

Set parameters $c=0.05$, $k=1$, $\alpha=1$, $A_1=0.3$, $A_2=1.5$, $\omega_1=1$ and $\omega_2=0.115$ . The whole period $T=2\pi/\mathrm{GCD}(\omega_1,\omega_2)=400\pi$ is many orders of magnitude larger than the period of the individual frequency component \cite{Chua1981}. As presented in Figure \ref{FIG:RHB}, the classical RHB method does not perform well in solving the multiple-frequency excitations case. We can conclude that the superposition of these two signals produces a periodic signal with a frequency of 0.005 (selected as the base frequency in the RHB method). So RHB method with 100 harmonics (RHB100) only considers the harmonic components of frequency 0.5 at most. Neglecting the higher frequency parts makes even the high-order RHB method shows a relatively big discrepancy from the benchmark result. Only when the harmonics cover up to frequency 1 (using RHB200 at least), it can realize a basic simulation for both low and high-frequency parts. To obtain a highly accurate solution, a high-order estimation form of the RHB method must be used, which is computationally expensive. Our newly purposed method tries to utilize fewer harmonics to realize more efficient and accurate estimation.

\begin{figure}
	\centering
	 \includegraphics[scale=0.5]{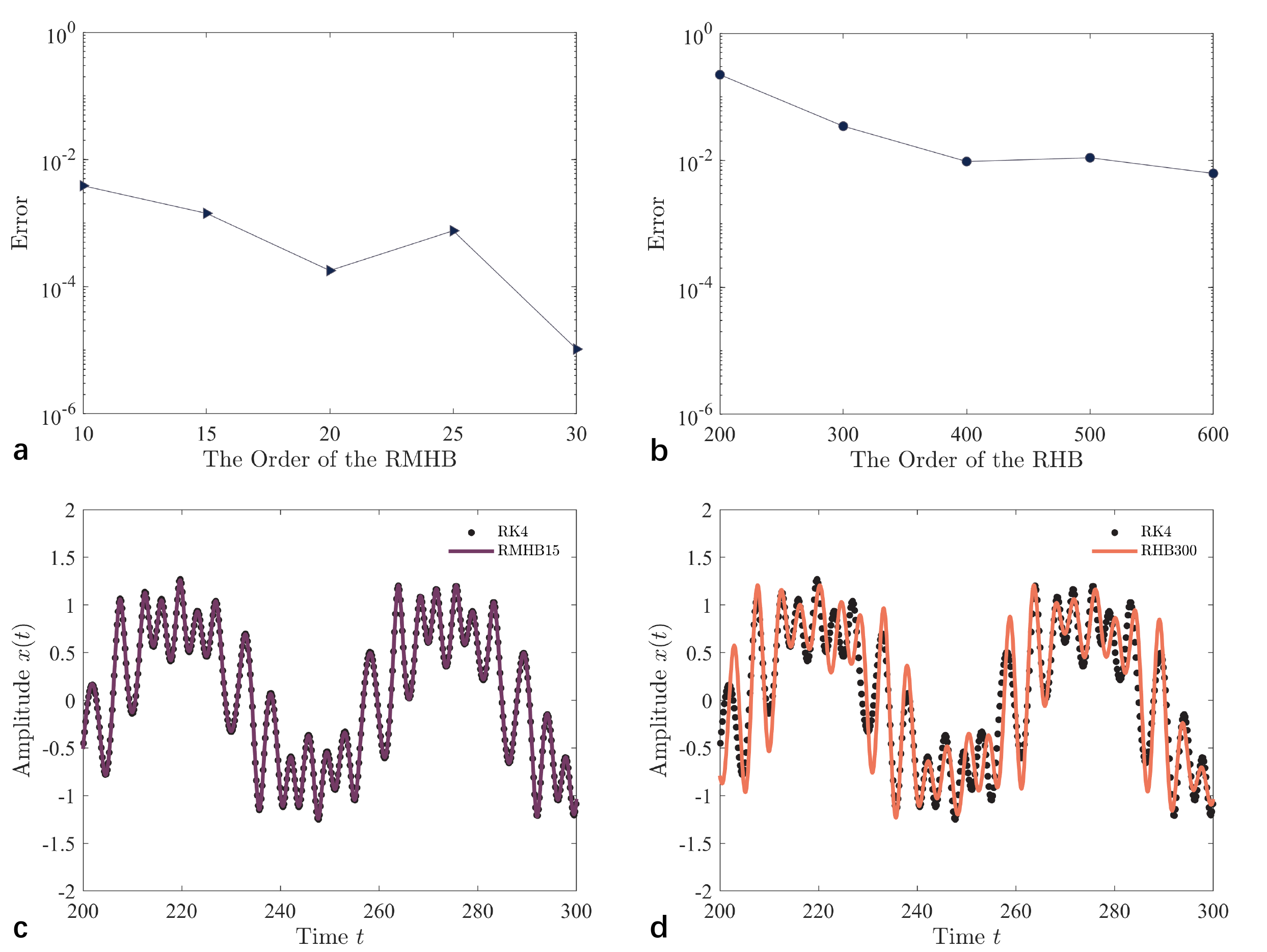}
	 \caption{Periodic analysis of Duffing equation $\ddot{x}+0.05\dot{x}+x+x^3=0.3\cos t+1.5\cos 0.115t$ using the RMHB and RHB methods. (a, b), Amplitude error curves of the RMHB method (a) and the RHB method (b) with different orders against the benchmark numerical result. (c, d), Time solution curves obtained by the RMHB15 (c) and the RHB300 (d).}
	 \label{FIG:RHBError}
\end{figure}

The initial value of the unknown variables of the RMHB method is assumed to be $\hat{x}_c(1,0)=0.2$, $\hat{x}_c(0,1)=0.8$. Figure \ref{FIG:RHBError}a shows that increasing the order can improve the computing accuracy of the RMHB method, the amplitude error of the RMHB30 is controlled at $10^{-5}$. While the RHB method has an insignificant effect on improving the accuracy. Table \ref{tab2} provides comparison results of the critical number of $M$, the amplitude error, and the computing time by using the RMHB and RHB methods. We find that the RMHB15 only accounts for half the computing time of the RHB300 but decreases the error by a factor of 25.

To sum up, for the strong nonlinear system with two external excitation inputs, the RHB method needs at least 200 orders to realize the amplitude estimation of the corresponding base frequencies $\omega_1$ and $\omega_2$. Due to the high-order estimation, the solving workload of the NAEs skyrockets, However, most harmonics are technically redundant. Therefore, the introduction of the RMHB method is to compress those unknown coefficients and discard redundant variables in the original RHB method to improve both the computational efficiency and accuracy as much as possible.

\begin{table}[htbp]
\centering
\footnotesize{
 \caption{\label{tab2}Results for various numerical methods when solving the multi-inputs Duffing equation}
 \begin{tabular}{cccc}
  \toprule
    Method & $M$ & Amplitude error & Computing time (s) \\ 
  \midrule
      RMHB1 & 801 & 0.1325 & 0.03 \\
    RMHB5 & 4001 & 0.0027 & 0.54\\
    RMHB15 & 12001 & 0.0014 & 50.29 \\
 \midrule
    RHB200 & 801 & 0.2230 & 9.00 \\
    RHB250 & 1001 & 0.2714 & 22.06 \\
    RHB300 & 1201 & 0.0344 & 124.28 \\
  \bottomrule
 \end{tabular}}
\end{table}
\subsection{Steady-State Periodic Aeroelastic Response Analysis of an Airfoil with  Multiple Base Frequencies}
\begin{figure}
	\centering
	 \includegraphics[scale=0.4]{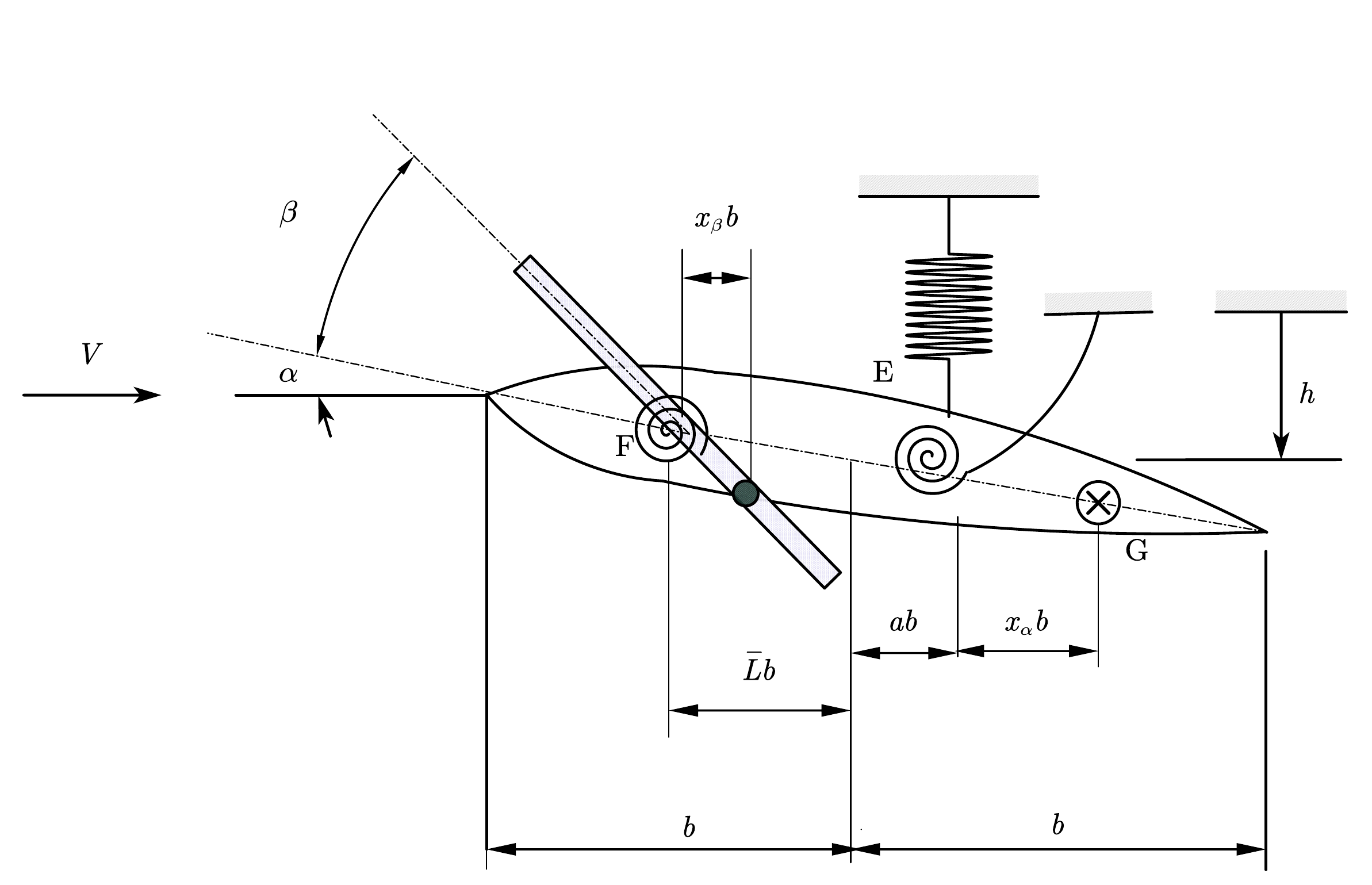}
	 \caption{Sketch of an airfoil with an external store}
	 \label{FIG:airfoil}
\end{figure}

A nonlinear aeroelastic system shown in Figure \ref{FIG:airfoil} is a two-dimensional airfoil with an external store. The airfoil section itself oscillates in both two directions of pitch and plunge. The plunge deflection is denoted by $h$, the pitch angle about the elastic axis is $\alpha$ and the varying pitch angle of the external store is $\beta$. Only the aerodynamic force acting on the airfoil is considered, while the effect on the external suspension is ignored. Taking the restoring forces into account as cubic nonlinearities, then the governing equations can be described as
\begin{equation}
    \mathbf{{M}\ddot{q}}+\mathbf{{C}\dot{q}}+\mathbf{Kq}+\mathbf{P}\left[ q_{1}^{3},q_{2}^{3},q_{3}^{3} \right] ^{\mathrm{T}}=0,
\end{equation}
the generalized coordinate vector is $\mathbf{q}=[ h, \alpha, \beta]^{\mathrm{T}}$, with $\mathbf{M}$, $\mathbf{C}$ and $\mathbf{K}$ are the mass, damping, and stiffness matrix respectively of size $3\times 3$. $\mathbf{P}$ is the parameter matrices distributing the nonlinear function. Those matrices are defined as 
\begin{equation}
    \mathbf{M}=\left[ \begin{matrix}
	\mu +\mu _{\beta}&		\mu x_{\alpha}+\mu _{\beta}x_{\beta}-\mu _{\beta}\bar{L}&		\mu _{\beta}x_{\beta}\\
	\mu x_a+\mu _{\beta}x_{\beta}-\mu _{\beta}\bar{L}&		\mu r_{\alpha}^{2}+\mu _{\beta}r_{\beta}^{2}+\mu _{\beta}\bar{L}^2-2\mu _{\beta}x_{\beta}\bar{L}&		\mu _{\beta}r_{\beta}^{2}-\mu _{\beta}x_{\beta}\bar{L}\\
	\mu _{\beta}x_{\beta}&		\mu _{\beta}r_{\beta}^{2}-\mu _{\beta}x_{\beta}\bar{L}&		\mu _{\beta}r_{\beta}^{2}\\
\end{matrix} \right],\nonumber 
\end{equation}

\begin{equation}
    \mathbf{C}=\left[ \begin{matrix}
	c_h&		0&		0\\
	0&		c_{\alpha}&		0\\
	0&		0&		c_{\beta}\\
\end{matrix} \right] ,\nonumber
\end{equation}

\begin{equation}
\mathbf{K}=\left[ \begin{matrix}
	\mu \left( \omega _h/\omega _{\alpha} \right) ^2&		2Q&		0\\
	0&		\mu r_{\alpha}^{2}-2\left( \bar{L}+a \right) Q&		0\\
	0&		0&		\mu _{\beta}r_{\beta}^{2}\left( \omega _{\beta}/\omega _{\alpha} \right) ^2\\
\end{matrix} \right] , \,\,
\mathbf{P}=\left[ \begin{matrix}
	k_{h3}&		0&		0\\
	0&		k_{\alpha 3}&		0\\
	0&		0&		k_{\beta 3}\\
\end{matrix} \right] ,
\nonumber
\end{equation}

with $Q=(V/b\omega_2)^2$ is a non-dimensional flow velocity, the physical meanings of other parameters in the equations above can be found in detail \cite{Liu2018,Chen2012}. The system parameters are chosen as $Q=8.0$, $\mu=12.8$, $\mu_{\beta}=4.0$, $x_{\alpha}=0.15$, $r_{\alpha}^{2}=0.3$, $r_{\beta}^{2}=0.89$, $\bar{L}=0.18$, $a=-0.41$, $b=0.118$, $c_h=0.2$, $c_\alpha=0.2$, $c_\beta=0$, $\omega_h=34.6$, $\omega_\alpha=88$ and $\omega_\beta=60$. Set nonlinear stiffness in the direction of the external store as $k_{h3}=k_{\alpha3}=0$ and $k_{\beta3}=10$. Through the fast Fourier transform (FFT) to analyze the amplitude spectrum, it is known that the steady-state periodic solution response of the autonomous system contains two base frequencies $f_1=0.0685$ and $f_2=0.0873$ \cite{Liu2018}. The RMHB algebraic system can be written as 
\begin{equation}
    \begin{aligned}
        \mathbf{M} \otimes \mathbf{I}_{2 p(p+1)+1}\left(\mathbf{E}^{*} \mathbf{E} \nabla^{2} \hat{\mathbf{q}}\right)+ & \mathbf{C} \otimes \mathbf{I}_{2 p(p+1)+1}\left(\mathbf{E}^{*} \mathbf{E} \nabla \hat{\mathbf{q}}\right)+ \\
        &\mathbf{K} \otimes \mathbf{I}_{2 p(p+1)+1}\left(\mathbf{E}^{*} \mathbf{E} \hat{\mathbf{q}}\right)+\mathbf{P} \otimes\mathbf{I}_{2 p(p+1)+1}(\mathbf{E}^{*} \tilde{\mathbf{f}}_M)=0,
    \end{aligned}
\end{equation}

with

\begin{equation}
    \nabla ^2=\mathbf{I}_N\otimes \left( \mathrm{diag}\left[ 0,\nabla _{1,0},\cdots ,\nabla _{m,n},\cdots ,\nabla _{0,p} \right] ^2 \right), \,\,
    \tilde{\mathbf{f}}_M=\left[ \begin{array}{c}
	\begin{array}{c}
	h^3\left( t_1 \right)\\
	\vdots\\
\end{array}\\
	\alpha^3\left( t_1 \right) \\
	\vdots\\
	\beta ^3\left( t_M \right) \\
\end{array} \right]  .\nonumber
\end{equation}

According to Theorem \ref{thm2}, the sampling time $T=1/\mathrm{GCD}(f_1,f_2)=10^4$, time domain nodes
\begin{equation}
    M>((\phi+1)p\cdot\mathrm{max}(f_1,f_2))/\mathrm{GCD}(f_1,f_2).\nonumber
\end{equation}

\begin{figure}
	\centering
	 \includegraphics[scale=0.3]{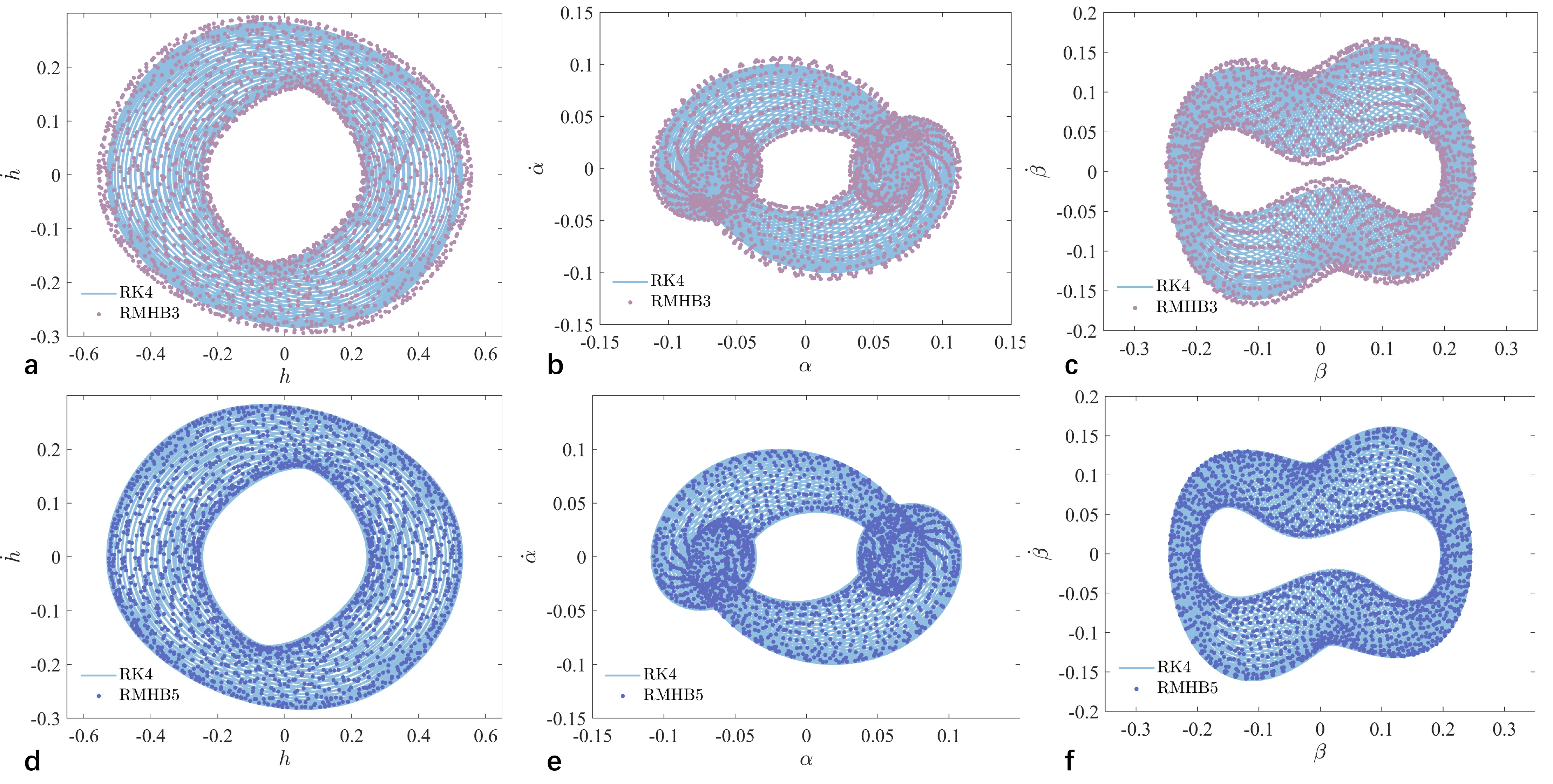}
	 \caption{Phase plot of quasi-periodic solutions obtained by the RMHB3 (a, b, c) and the RMHB5 (d, e, f).}
	 \label{FIG:RMHBQP}
\end{figure}

Figure \ref{FIG:RMHBQP} exhibits the phase plot obtained by the RMHB3 (see Figure \ref{FIG:RMHBQP}a-\ref{FIG:RMHBQP}c) and RMHB5 (see Figure \ref{FIG:RMHBQP}d-\ref{FIG:RMHBQP}f) respectively, benchmarked with the RK4 method. Engineering problems are generally multi-dimensional, which limits the HB practical use in complex dynamics studies. But the RMHB method builds the system of  NAEs more efficiently, through matrix operations rather than symbolic operations. Table \ref{tab3} shows the computing error and efficiency for different orders of the RMHB method. Here we denote $|\Delta A_h|$, $|\Delta A_\alpha|$ and $|\Delta A_\beta|$ are their amplitude errors corresponding to $h$, $\alpha$ and $\beta$. Even when solving the steady-state response of a multi-DOF system with high-order estimation, its efficiency and accuracy can be ensured. Furthermore, the computation precision can be improved by taking more harmonics into account. Nevertheless, it can be also found that the difference between the amplitude errors corresponding to $p=5$ and $p=7$ is not apparent, in other words,  when the number of reserved harmonics is enough, the improvement effect of accuracy is insignificant.

\begin{table}[htbp]
\centering
\footnotesize{
 \caption{\label{tab3}Comparison of amplitude error and efficiency of the RMHB method with different orders}
 \begin{tabular}{ccccc}
  \toprule
    Order & $\left| \Delta A_h \right|$ & $\left| \Delta A_\alpha \right|$ & $\left| \Delta A_\beta \right|$ & Calculation time (s) \\ 
  \midrule
    1 & 0.0992 & 0.0198 & 0.0340 & 3.24 \\
    3 & 0.0266 & 0.0040 & 0.0032 & 37.60 \\
    5 & 0.0060 & 0.0016 & $6.65\times 10^{-4}$ & 150.95\\
    7 & 0.0035 & 0.0011 & $2.80\times 10^{-4}$ & 464.28\\
  \bottomrule
 \end{tabular}}
\end{table}

\section{Conclusion}
In this study, the reconstructive harmonic balance (RMHB) method is introduced to obtain quasi-periodic responses of a nonlinear dynamic system robustly and efficiently, which extends the RHB to multiple base harmonics cases. In this method, the resultant NAEs can be simply derived, thus avoiding complicated symbolic operations. We also reveal the theoretical mechanism of the aliasing phenomenon in the multiple harmonic balance computations and thus put forward sampling rules to eliminate aliasing for different multiple base frequencies problems.

First, the forced Van der Pol equation is computed to elucidate how the present method tackles the quasi-periodic responses. If the ratio of base frequencies is irrational, the dynamic system will generate a quasi-periodic response. We demonstrate that the RMHB method works more efficiently and ensures convergence, which archives at least 1000 times speed up fast in the computing times than the HDHB method (with two base frequencies). Furthermore, the changes of non-zero elements and the maximum element value in the aliasing matrix are checked to investigate the effect of de-aliasing. Studies have shown that the RMHB method can minimize the aliasing error (error between the RMHB and MHB method) with the increase in the $T$ and $M$. The aliasing error can be reduced to $10^{-5}$ with $T=M=5\times 10^5$. With effectiveness and high precision, the RMHB method with multiple base frequencies could be applicable to other nonlinear dynamical systems, especially in detecting and analyzing quasi-periodic solutions.

Second, Duffing oscillators with two input signals are explored. The input frequencies differ by many orders of magnitude but the frequency ratio is rational. Through the analysis of the distribution of solutions, we find that the aliasing can be fully eliminated by choosing the proper sampling period $T$ and the number of collocations $M$. Whereas it will produce additional non-physical solutions if the conditional equivalence no longer holds. Besides the RMHB method improves efficiency by compressing redundant variables. In fact, there are only several harmonics taking the main part of solutions in the classical HB-like methods. Thus, it allows fewer frequency domain variables to achieve more accurate solutions. Being an optimal reconstruction of the MHB method, we find that the computational error of the RMHB method can be controlled to the order of $10^{-5}$, which is $10^3$ more accurate than the RHB. Third, the nonlinear response of an airfoil with an external store is analyzed. The quasi-periodic solutions obtained by the present method are in excellent consistency with the results provided by numerical integration. Besides the credible accuracy, coupled differential equations with strong nonlinearity can be effectively handled by the RMHB method. But the use of the MHB method is severely limited by symbolic operations when computing such problems.

The computational performance of the RMHB method for the quasi-periodic response problem of nonlinear systems has advantages over existing methods, which is expected to become a fast and accurate method for solving multiple base frequencies with the theoretical meaning of de-aliasing.  The semi-analytical solution obtained by the RMHB method provides us with a convenient way to study the properties of the steady-state responses across multiple disciplines. However, the present method has a limitation: since both the collocation matrix $\mathbf{E}$ and the transformation matrix $\mathbf{E}^*$ are explicitly defined, this newly proposed approach requires prior information about those base frequencies. But the frequency components are often undisclosed for autonomous systems, thus more dedicated improvements \cite{Zheng2022} will be considered in our future work. 

\appendix
\setcounter{equation}{0}
\renewcommand{\theequation}{A.\arabic{equation}}
\section{Cubic Nonlinear Expansion}
For truncation order $p=1$, the approximated function is 
\begin{equation}
    x=\hat{x}_c(0,0) +\hat{x}_c(1,0)\cos \omega _1t+\hat{x}_s(1,0)\sin \omega_1t
       +\hat{x}_c(0,1)\cos \omega _2t+\hat{x}_s(0,1)\sin \omega _2t. 
\end{equation}
The Fourier coefficient vector for cubic nonlinear terms is $\mathbf{\hat{f}}(\mathbf{\hat{x}})$, which is manually sorted as
\begin{equation}
    \hat{\mathbf{f}}=[\hat{f}_{c}(0,0)\,\,\hat{f}_{c}(1,0)\,\,\hat{f}_{s}(1,0)\,\,\hat{f}_{c}(0,1)\,\, \hat{f}_{s}(0,1)]^{\mathrm{T}},
\end{equation}
where
\begin{align*}
    \hat{f}_{c}(0,0)=\hat{x}_{c}(0,0)^{3} +\frac{3}{2}\hat{x}_{c}(0,0)\hat{x}_{c}(1,0)^{2}
    &+\frac{3}{2}\hat{x}_{c}(0,0)\hat{x}_{s}(1,0)^{2} \\
    &+\frac{3}{2}\hat{x}_{c}(0,0)\hat{x}_{c}(0,1)^{2}+\frac{3}{2}\hat{x}_{c}(0,0)\hat{x}_{s}(0,1)^{2},
\end{align*}

\begin{align*}
    &\hat{f}_{c}(1,0)= 3\hat{x}_{c}(1,0)\hat{x}_{c}(0,0)^{2}+\frac{3}{4}\hat{x}_{c}(1,0)^{3}+\frac{3}{2}\hat{x}_{c}(1,0)\hat{x}_{c}(0,1)^{2}+
    \frac{3}{4}\hat{x}_{c}(1,0)\hat{x}_{s}(1,0)^{2},
\end{align*}

\begin{align*}
    \hat{f}_{s}(1,0)= 3\hat{x}_{s}(1,0)\hat{x}_{c}(0,0)^{2}+ \frac{3}{4}\hat{x}_{s}(1,0) \hat{x}_{c}(1,0)^{2} & +\frac{3}{4}\hat{x}_{s}(1,0)^{3} \\ 
    & + \frac{3}{2}\hat{x}_{s}(1,0)\hat{x}_{c}(0,1)^{2}+\frac{3}{2}\hat{x}_{s}(1,0)\hat{x}_{s}(0,1)^{2},
\end{align*}

\begin{align*}
    \hat{f}_{c}(0,1)= 3\hat{x}_{c}(0,1)\hat{x}_{c}(0,0)^{2}+\frac{3}{2}\hat{x}_{c}(0,1) \hat{x}_{c}(1,0)^{2} &+\frac{3}{2}\hat{x}_{c}(0,1)\hat{x}_{s}(1,0)^{2} \\
    &+\frac{3}{4}\hat{x}_{c}(0,1)^{3}+\frac{3}{4}\hat{x}_{c}(0,1)\hat{x}_{s}(0,1)^{2},
\end{align*}

\begin{align*}
    \hat{f}_{s}(0,1)= 3\hat{x}_{s}(0,1)\hat{x}_{c}(0,0)^{2}+\frac{3}{2}\hat{x}_{s}(0,1) \hat{x}_{1,0}^{2} &+\frac{3}{2}\hat{x}_{s}(0,1) \hat{x}_{s}(1,0)^{2} \\
    &+\frac{3}{4}\hat{x}_{s}(0,1)\hat{x}_{c}(0,1)^{2}+\frac{3}{4}\hat{x}_{s}(0,1)^{3}.
\end{align*}

The corresponding Fourier coefficients of higher-order harmonics obtained by cubic expansion are:
\begin{align*}
    &\hat{f}_{c}(2,0)=\frac{3}{2}\hat{x}_{c}(0,0)\hat{x}_{c}(1,0)^{2}-\frac{3}{2}\hat{x}_{c}(0,0) \hat{x}_{s}(1,0)^{2},\quad
    \hat{f}_{s}(2,0)= 3\hat{x}_{c}(0,0)\hat{x}_{c}(1,0)\hat{x}_{s}(1,0),\\
    &\hat{f}_{c}(3,0)= \frac{1}{4}\hat{x}_{c}(1,0)^{3}-\frac{3}{4}\hat{x}_{c}(1,0) \hat{x}_{s}(1,0)^{2},\quad
    \hat{f}_{s}(3,0)=\frac{3}{4}\hat{x}_{s}(1,0)\hat{x}_{c}(1,0)^{2}-\frac{1}{4}\hat{x}_{s}(1,0)^{3},\\
    &\hat{f}_{c}(-2,1)=\frac{3}{4}\hat{x}_{c}(1,0)\hat{x}_{c}(0,1)^{2}+\frac{3}{2}\hat{x}_{c}(0,1) \hat{x}_{s}(0,1)\hat{x}_{s}(1,0),\\
    &\hat{f}_{s}(-2,1)=\frac{3}{4}\hat{x}_{s}(1,0)\hat{x}_{c}(0,1)^{2}-\frac{3}{2}\hat{x}_{c}(1,0) \hat{x}_{c}(0,1)\hat{x}_{s}(0,1)-
     \frac{3}{4}\hat{x}_{s}(1,0)\hat{x}_{s}(0,1)^{2},\\
    &\hat{f}_{c}(-1,1)= 3\hat{x}_{c}(0,0)\hat{x}_{c}(1,0) \hat{x}_{c}(0,1)+3\hat{x}_{c}(0,0) \hat{x}_{s}(1,0)\hat{x}_{s}(0,1),\\
    &\hat{f}_{s}(-1,1)=-3\hat{x}_{c}(0,0)\hat{x}_{c}(1,0)\hat{x}_{s}(0,1)+3\hat{x}_{c}(0,0) \hat{x}_{s}(1,0)\hat{x}_{c}(0,1),\\
    &\hat{f}_{c}(1,1)=3\hat{x}_{c}(0,0) \hat{x}_{c}(1,0)\hat{x}_{c}(0,1)-3\hat{x}_{c}(0,0) \hat{x}_{s}(1,0) \hat{x}_{s}(0,1),\\
    &\hat{f}_{s}(1,1)=3\hat{x}_{c}(0,0)\hat{x}_{c}(1,0)\hat{x}_{s}(0,1)+3\hat{x}_{c}(0,0) \hat{x}_{s}(1,0) \hat{x}_{c}(0,1),\\
    &\hat{f}_{c}(2,1)=\frac{3}{4}\hat{x}_{c}(0,1)\hat{x}_{c}(1,0)^{2}-\frac{3}{2}\hat{x}_{c}(1,0)\hat{x}_{s}(1,0)\hat{x}_{s}(0,1)-
    \frac{3}{4}\hat{x}_{c}(0,1)\hat{x}_{s}(1,0)^{2}, \\
    &\hat{f}_{s}(2,1)=\frac{3}{4} \hat{x}_{s}(0,1)\hat{x}_{c}(1,0)^{2}+\frac{3}{2} \hat{x}_{c}(1,0) \hat{x}_{s}(1,0) \hat{x}_{c}(0,1)-
    \frac{3}{4}\hat{x}_{s}(0,1) \hat{x}_{s}(1,0)^{2},\\
    &\hat{f}_{c}(-1,2)=\frac{3}{4}\hat{x}_{c}(1,0) \hat{x}_{c}(0,1)^{2} +\frac{3}{2} \hat{x}_{s}(1,0) \hat{x}_{c}(0,1)\hat{x}_{s}(0,1),\\
    &\hat{f}_{s}(-1,2)=\frac{3}{2} \hat{x}_{c}(1,0)\hat{x}_{c}(0,1)\hat{x}_{s}(0,1)-\frac{3}{4} \hat{x}_{s}(1,0)\hat{x}_{c}(0,1)^{2} +
    \frac{3}{4} \hat{x}_{s}(1,0)\hat{x}_{s}(0,1)^{2},\\
    &\hat{f}_{c}(0,2)=\frac{3}{2}\hat{x}_{c}(0,0)\hat{x}_{c}(0,1)^{2}-\frac{3}{2}\hat{x}_{c}(0,0)\hat{x}_{s}(0,1)^{2},\quad 
    \hat{f}_{s}(0,2) = 3\hat{x}_{c}(0,0)\hat{x}_{c}(0,1)\hat{x}_{s}(0,1),\\
    &\hat{f}_c(1,2) =\frac{3}{4}\hat{x}_c(1,0)\hat{x}_c(0,1) ^2 -\frac{3}{2} \hat{x}_s( 1,0 ) \hat{x}_c(0,1) \hat{x}_s(0,1)-
    \frac{3}{4}\hat{x}_{c}(1,0) \hat{x}_{s}(0,1)^{2},\\
    &\hat{f}_{s}(1,2)=\frac{3}{2}\hat{x}_{c}(1,0)\hat{x}_{c}(0,1)\hat{x}_{s}(0,1)+\frac{3}{4} \hat{x}_{s}(1,0) \hat{x}_{c}(0,1)^{2}-
    \frac{3}{4} \hat{x}_{s}(1,0) \hat{x}_{s}(0,1)^{2},\\
    &\hat{f}_{c}(0,3) = \frac{1}{4}\hat{x}_{c}(1,0)^{3}-\frac{3}{4}\hat{x}_{c}(1,0)\hat{x}_{s}(1,0)^{2},
    \quad \hat{f}_{s}(0,3) = \frac{3}{4}\hat{x}_{s}(1,0) \hat{x}_{c}(1,0)^{2}-\frac{1}{4} \hat{x}_{s}(1,0)^{3}.
\end{align*}

\newpage
\bibliographystyle{unsrt}
\bibliography{ref}

\end{document}